\documentclass[11pt]{article}
\usepackage{epsfig,epstopdf,graphicx}
\usepackage{amsmath}
\usepackage{amsfonts}
\usepackage{amssymb}
\usepackage[svgnames]{xcolor}
\usepackage{bm,graphicx}
\usepackage{caption}

\newtheorem{thm}{Theorem}[section]
\newtheorem{lem}[thm]{Lemma}
\newtheorem{prop}[thm]{Proposition}

\newtheorem{defn}[thm]{Definition}
\newtheorem{example}[thm]{Example}

\newtheorem{rem}[thm]{Remark}

\numberwithin{equation}{section}
\newcommand{\R}{\mathbb R}

\newcommand{\rr}{\mathbb R}
\newcommand{\E}{\mathbb E}

\newcommand{\D}{\mathbb D}
\newcommand{\I}{\mathbb I}

\newcommand{\N}{\mathbb N}

\newcommand{\ip}[2]{{\langle #1,#2\rangle}}
\newcommand{\Var}{{\operatorname{Var}}}
\renewcommand{\Im}{{\mathcal I}}
\newcommand{\1}{\mathbf 1}

\begin{document}
%\sloppy
\title%[Tempered Fractionally Integrated Processes]
%{How does tempering affect sample path properties and stochastic integration with respect to fractional L\'evy processes?} %Tempered fractional L\'{e}vy processes }
{On fractional L\'evy processes: tempering, sample path properties and stochastic integration}

\author{B.\ Cooper Boniece \\ Department of Mathematics and Statistics \\ Washington University in St.\ Louis \and Gustavo Didier \\ Mathematics Department\\ Tulane University \and Farzad Sabzikar\\ Department of Statistics\\ Iowa State University}
\date{\today \\  \small
\vskip.2cm
}

\maketitle

\begin{abstract}
We define two new classes of stochastic processes, called tempered fractional L\'{e}vy process of the first and second kinds (TFLP and TFLP \textit{II}, respectively). TFLP and TFLP \textit{II} make up very broad finite-variance, generally non-Gaussian families of transient anomalous diffusion models that are constructed by exponentially tempering the power law kernel in the moving average representation of a fractional L\'{e}vy process. Accordingly, the increment processes of TFLP and TFLP \textit{II} display semi-long range dependence. We establish the sample path properties of TFLP and TFLP \textit{II}. We further use a flexible framework of tempered fractional derivatives and integrals to develop the theory of stochastic integration with respect to TFLP and TFLP \textit{II}, which may not be semimartingales depending on the value of the memory parameter and choice of marginal distribution. %A statistical study of geophysical flow data shows that a second order, non-Gaussian tempered model provides a good fit.
\end{abstract}

%This paper introduces new infinitely divisible stochastic L\'{e}vy processes based on tempered fractional calculus which we call them tempered fractional L\'{e}vy processes of the first and second kind. These processes have semi-log memory. That is the covariance function resemble the covariance function of a long memory model for some moderate lags but eventually decays exponentially fast. Some essential properties including the covariance structure, sample paths behavior, and stochastic integrations with respect to the tempered fractional processes are presented. Finally, the application of these models for real data will be discussed.

\section{Introduction}

%\GDcomment{Interesting issues: $(i)$ (non)equivalence of moving average and harmonizable representations?; $(ii)$ in general, what are the natural conditions for the equivalence of the multiple definitions of semi-LRD?}

In this paper, we define two new classes of stochastic processes, called \textit{tempered fractional L\'{e}vy processes} of the first and second kinds (TFLP and TFLP \textit{II}, respectively). TFLP and TFLP \textit{II} make up very broad finite-variance, generally non-Gaussian transient anomalous diffusion models, i.e., their second order properties qualitatively change over time. They are constructed by exponentially tempering the power law kernel in the moving average representation of a fractional L\'{e}vy process (FLP). In particular, their increment processes exhibit \textit{semi-long range dependence} (semi-LRD) in the sense of \cite{giraitis}, namely, their autocovariance functions decay hyperbolically over small lags and exponentially fast over large lags (see \eqref{e:semi-LRD}). We establish the sample path regularity of TFLPs. Turning to stochastic analysis, we use a flexible framework of tempered fractional derivatives and integrals to develop the theory of stochastic integration with respect to TFLP and TFLP \textit{II}, which may not be semimartingales depending on the value of the memory parameter and choice of marginal distribution. %A statistical study of geophysical flow data demonstrates the need for non-Gaussian transiently anomalous physical models.

Fractional, or non-Markovian, stochastic processes naturally emerge in many fields of science, technology and engineering (see, e.g., \cite{ mandelbrot:vanness:1968,ciuciu:abry:he:2014,Foufoula94,ivanov1999,Mandelbrot1974,taqqu97}). They provide the mathematical framework for what is called \textit{scale-free analysis} \cite{mandelbrot:vanness:1968,flandrin:1992,wornell:oppenheim:1992}. Rather than focusing on the detection of a small number of characteristic scales, in scale-free analysis it is assumed that the phenomenological dynamics are driven by a large continuum of time scales usually related by means of a \textit{power law}. %Self-similar processes are among the most popular fractional or scale invariant models~\cite{Samorodnitsky1994}. A stochastic processes is called \textit{self-similar} when time dilation is probabilistically neutral up to a power law, namely, $\{X(t)\}_{ t \in \R}  \stackrel{\textnormal{f.d.d.}}{=} \{c^H X(t/c)\}_{t \in \R}$, $c > 0$, where $H > 0$ is called the Hurst exponent and $\stackrel{\textnormal{f.d.d.}}{=} $ denotes the equality of finite-dimensional distributions.
A cornerstone class of scale invariant processes is fractional Brownian motion (FBM), i.e., the only Gaussian, self-similar, stationary increment process \cite{embrechts:maejima:2002,pipiras:taqqu:2017}. The literature on fractional processes is now voluminous; see, e.g., \cite{beran:feng:ghosh:kulik:2013,dobrushin:major:1979,granger:joyeux:1980,moulines:roueff:taqqu:2008,taqqu:1975,taqqu:1979,SamorodnitskyTaqqu,bardet:tudor:2014,clausel:roueff:taqqu:tudor:2014:waveletestimation}.

\begin{figure}[h!]%
\begin{center}
\vskip-10pt
\includegraphics[scale=0.5]{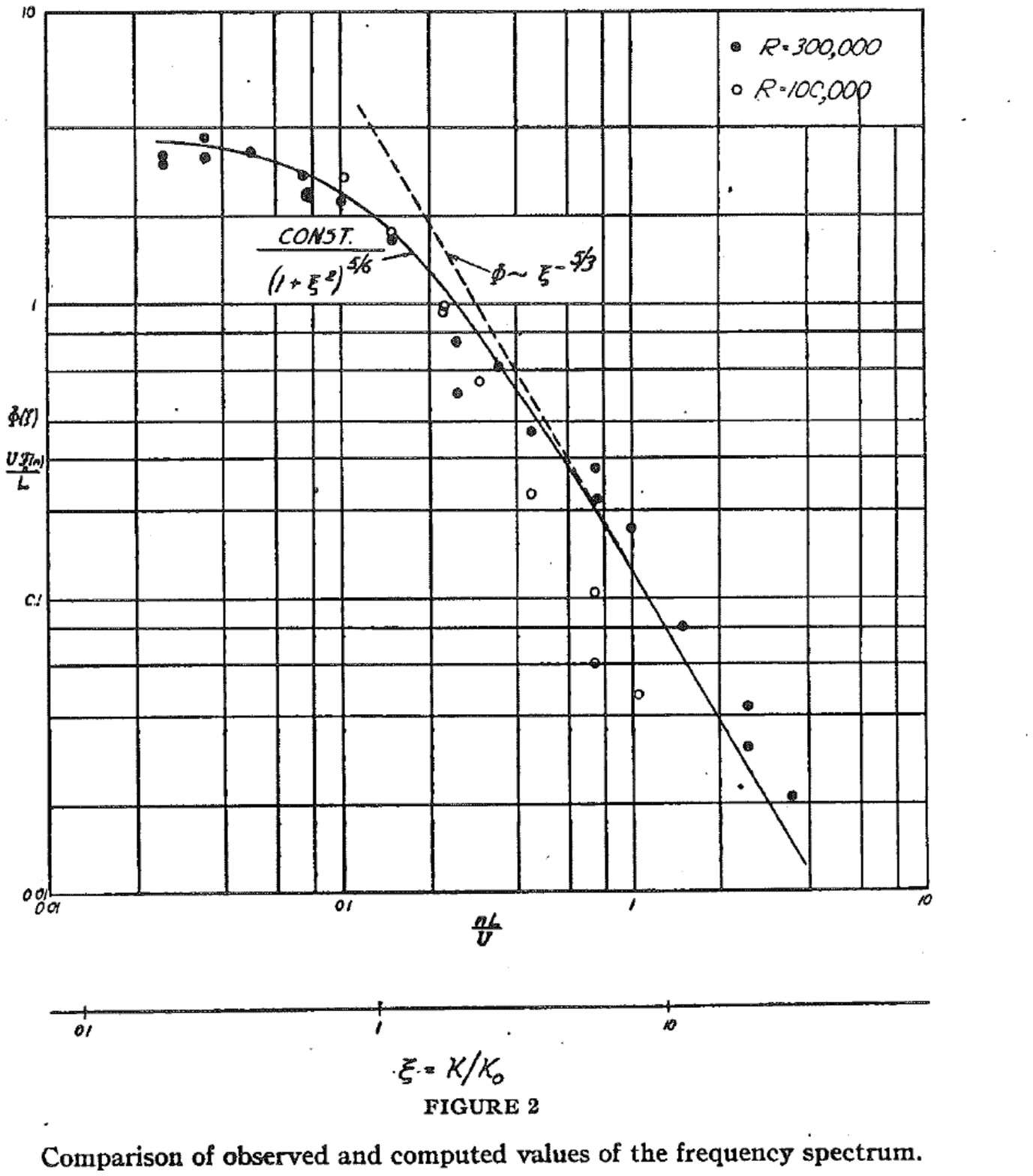}
\end{center}
\label{fig:vonkarman}
\caption{\textbf{The Von K\'{a}rm\'{a}n spectral density (curved line) versus Kolmogorov's 5/3 law (straight line).} Kolmogorov's classical theory \cite{kolmogorov:1940,Kolmogorovturbulence,friedlander:topper:1961,shiryaev:1999} posits that the energy spectrum in the inertial range is universal and given by the frequency domain power law $\omega^{-5/3}$. However, in the production (low frequency) range, turbulence is not universal, which may lead to transient behavior. In the Von K\'{a}rm\'{a}n model of continuous wind gusts~\cite{vonkarman1,DOD:2004,penner:williams:libby:nemat-nasser:2009,Beaupuits,Jang,Norton}, the framework favored by the U.S.\ Department of Defense in aircraft design, the spectral density is proportional to $(\lambda^2 +\omega^2)^{-5/6}$ with $\lambda = 1$ (see \cite{davenport:1961,norton:wolff:1981,li:kareem:1990,beaupuits:etal:2004,DOD:2004}). The tempering parameter $\lambda$ dampens down, in the low frequency limit, the power law behavior universally valid for the inertial range. The spectral density of tempered fractional L\'{e}vy noise $I\!I$, the increment process of TFLP $I\!I$, is of the Von K\'{a}rm\'{a}n type (see Proposition \ref{p:TFLNII}).}\label{fig:vonkarman}
\end{figure}
In many empirical settings, power law behavior is expected to hold only within a range of scales, out of which the observed dynamics qualitatively change, possibly to different power law behavior or simply non-fractional stationarity. In anomalous diffusion modeling, this is typically reflected in the behavior of the so-named mean squared displacement (MSD)
\begin{equation}\label{e:MSD}
\E X^2(t) \approx C t^{\vartheta}, \quad C,\vartheta \geq 0,
\end{equation}
of the particle position $X(t)$ over a time interval $T \ni t$, where the instances $\vartheta = 1$ and $\vartheta \neq 1$ correspond to classical and anomalous behavior, respectively (e.g., \cite{metzler:klafter:2000,kou:2008,sokolov:2008,didier:mckinley:hill:fricks:2012,grebenkov:vahabi:bertseva:forro:jeney:2013,zhang:crizer:schoenfisch:hill:didier:2018}). In the physics literature, a particle is said to undergo \textit{transient anomalous diffusion} when the value of the exponent $\vartheta$ in \eqref{e:MSD} changes over different time intervals (e.g., \cite{piryatinska:sanchev:woyczynski:2005,stanislavsky:weron:weron:2008,baeumer:meerschaert:2010,sandev:chechkin:kantz:metzler:2015,wu:deng:barkai:2016,chen:wang:deng:2017,liemert:sandev:kantz:2017,chen:wang:deng:2018}). Transience may appear in several contexts such as in nanobiophysics~\cite{saxton:2007,molina:sandev:safdari:pagnini:chechkin:metzler:2018} and particle dispersion \cite{taylor:1922,xia:francois:punzmann:shats:2013}. It also arises as a consequence of accounting for the energy spectrum of turbulence in the low frequency range, leading to the so-named Davenport-- \cite{boniece:didier:sabzikar:2018} or Von K\'{a}rm\'{a}n--type spectra (see Figure \ref{fig:vonkarman}).

\begin{figure}[h!]
    \centering
    \begin{minipage}{.5\linewidth}
        \centering
        \includegraphics[width=1\linewidth]{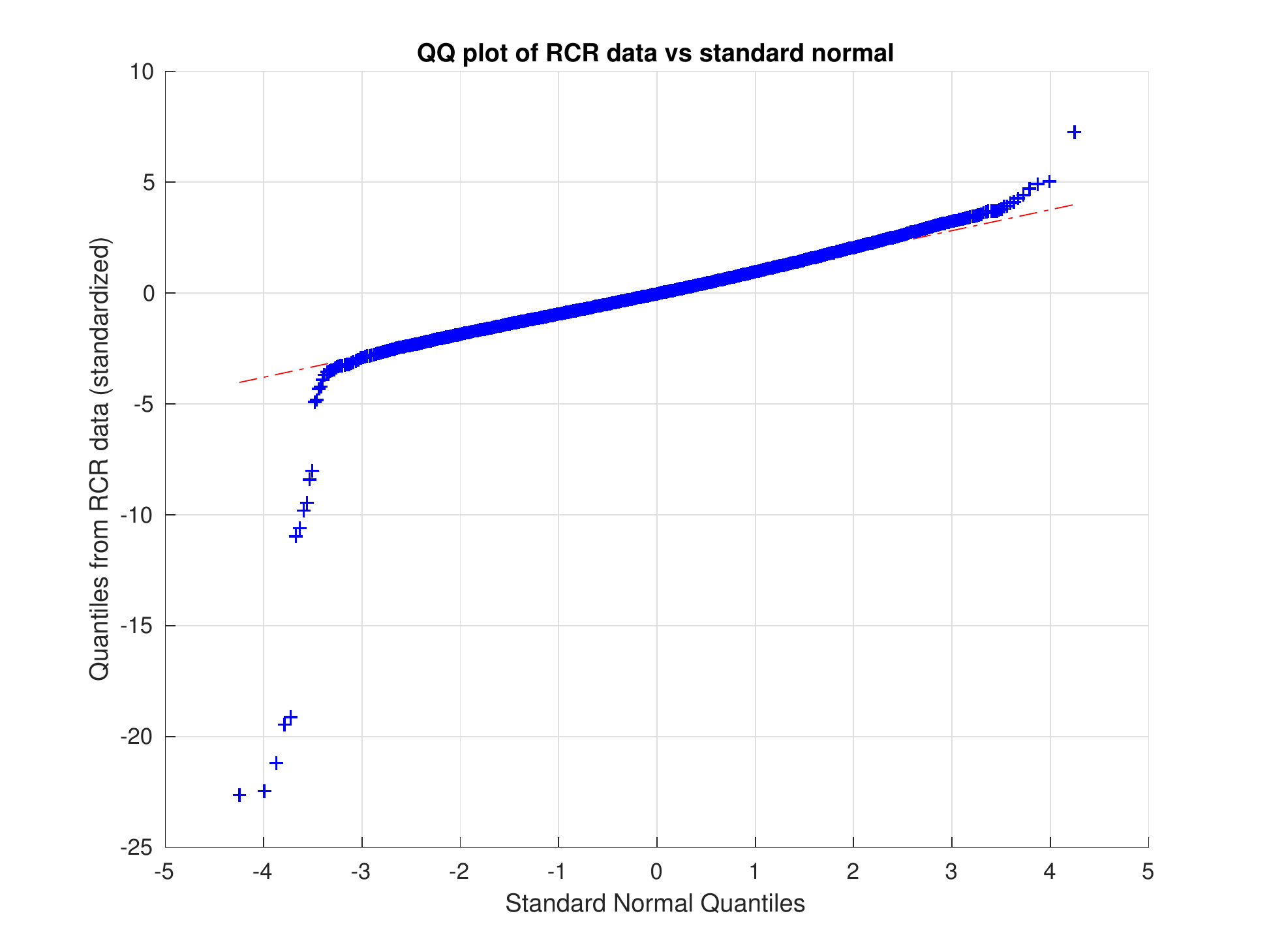}
      %  \caption{}

    \end{minipage}%
    %\begin{minipage}{.5\textwidth}
%        \centering
%        \includegraphics[width=.7\linewidth]{simulated_QQ1_closeup}
%%        \caption{}
%        \label{fig:prob1_6_1}
%    \end{minipage}
    \caption{\label{fig:RCR1} \textbf{Non-Gaussianity in river flow turbulence}.  Data on turbulent supercritical flow in the Red Cedar River, a fourth-order stream in Michigan, USA, was collected and kindly provided by Prof.\ Mantha S.\ Phanikumar, from Michigan State University. The measurements ($n=46080$ points) were made at a sampling rate of 50 Hz using a 16 MHz Sontek Micro-ADV (Acoustic Doppler Velocimeter) on May 26, 2014. The data is modeled in \cite{meerschaert:sabzikar:phanikumar:zeleke:2014} in the Fourier and in \cite{boniece:sabzikar:didier:2018,boniece:didier:sabzikar:2018} in the wavelet domains. The qq--plot, shown above, further reveals the conspicuous non-Gaussianity of the sample tails.
     }
       \end{figure}

Tempered FBM of the first and second kinds (TFBM \cite{Meerschaertsabzikar} and TFBM $I\!I$ \cite{SurgailisFarzadTFSMII}, respectively) are transient anomalous diffusion models. For TFBM, the MSD in \eqref{e:MSD} goes from $\vartheta > 0$ over small time scales to $\vartheta = 0$ over large scales, as in geophysical flows \cite{meerschaert:zhang:baeumer:2008,meerschaert:sabzikar:phanikumar:zeleke:2014}. By contrast, for TFBM $I\!I$, it shifts from anomalous over small scales to regular ($\vartheta = 1$) over large scales, as in viscoelastic diffusion (cf.\ \cite{fricks:yao:elston:forest:2009,francois:xia:punzmann:combriat:shatz:2015,xia:francois:punzmann:shats:2014}). Accordingly, the autocovariance functions $\gamma$ of the increments of both TFBM and TFBM $I\!I$ have the related property of semi-LRD, i.e.,
\begin{equation}\label{e:semi-LRD}
\gamma(h) \sim C \frac{|h|^{\delta}}{e^{\lambda |h|}}, \quad \lambda > 0, \quad \delta > -\frac{3}{2}, \quad |h| \rightarrow \infty,
\end{equation}
where $\lambda > 0$ is called the tempering parameter (see also Remark \ref{TFBMBSSdifference1} on the related literature on L\'evy semistationary processes). Moreover, like FBM vis-\`{a}-vis the Kolmogorov spectrum in the inertial range, TFBM $I\!I$ \cite{Meerschaertsabzikar,MeerschaertsabzikarSPA} is a Gaussian model that displays a von K\'{a}rm\'{a}n--type spectrum. Due to their appeal in applications, TFBMs have recently attracted considerable research efforts~\cite{zeng:yang:chen:2016,chen:wang:deng:2017}. In \cite{boniece:didier:sabzikar:2018,boniece:sabzikar:didier:2018}, wavelets are used in the construction of the first statistical method for TFBM as a model of geophysical flow turbulence. Nevertheless, there is abundant phenomenological evidence of non-Gaussian behavior, especially in terms of tail distributions. This is true, for example, for the velocity and velocity derivative processes in wind turbulence \cite{barndorff-nielsen:1977,barndorff-nielsen:1979,barndorff-nielsen:jensen:sorensen:1989,barndorff-nielsen:jensen:sorensen:1990,barndorff-nielsen:jensen:sorensen:1993,skyum:christiansen:blaesild:1996} or returns to financial assets \cite{barndorff-nielsen:1997}; see also Figure \ref{fig:RCR1}.  %\BCBcomment{In truth I don't think adding the QQ-plot for the RCR data adds any value since we already mention the existence of non-Gaussian features of turbulence here.}.
Accordingly, many authors have developed several other classes of tempered non-Gaussian stochastic processes such as tempered fractional stable or tempered Hermite processes~\cite{SurgailisFarzadTFSMII,sabzikar:2015}, and tempered stable processes \cite{rosinski:2007,baeumer:meerschaert:2010,bianchi:rachev:kim:fabozzi:2010,gajda:magdziarz:2010,rosinski:sinclair:2010,kawai:masuda:2012,kuchler:tappe:2013}.

The family of fractional L\'{e}vy processes (e.g., \cite{benassi:cohen:istas:2002,brockwell:marquart:2005,Marquardt,lacaux:loubes:2007,bender:marquardt:2008}) has become popular in physical modeling since it provides a second order non-Gaussian framework displaying fractional covariance structure~\cite{barndorff-nielsen:schmiegel:2008,suciu:2010,magdziarz:weron:2011,zhang:li:zhang:2015,xu:li:zhang:li:kurths:2016}. In this paper, we construct the classes of TFLP and TFLP $I\!I$, which are families of tempered fractional processes with finite-variance, infinitely divisible finite-dimensional distributions. While FLP (including FBM) is only well-defined for memory parameter values $d \in (-1/2,1/2)$~\cite{Fink,Marquardt}, TFLPs are well-defined for every $d>-1/2$ due to the tempering effect of the exponential function in their kernels. We establish their second order and sample path regularity properties (see Propositions \ref{prop:covariance_TFLP}, \ref{prop:TFLN}, \ref{prop:THPcovariance} and \ref{p:TFLNII} and Theorems \ref{t:TFLP_sample_path} and \ref{t:TFLPI bounds}). In our analysis, continuous modifications of TFLP and TFLP \textit{II} can also be obtained, under conditions, by means of improper Riemann integral representations (Propositions \ref{p:modification_Theorem_TFLP1} and \ref{modification_Theorem_TFLP2}; see also Bender et al.\ \cite{bender:knobloch:oberacker:2015} for related results in a general martingale-driven framework). In particular, our results show that TFLP and TFLP $I\!I$ can be viewed as non-Gaussian transient anomalous diffusion models whose second order properties generalize those of TFBM and TFBM $I\!I$, respectively (see also Example \ref{ex:simulation} and Figures \ref{fig:TFLPI_paths}, \ref{fig:TFLPII_paths} on the effect of non-Gaussian noise distributions on sample path behavior).

Physical models of transient phenomena are often based on Langevin-type stochastic differential equations; see, for example, \cite{molina:sandev:safdari:pagnini:chechkin:metzler:2018} on the transient MSD of solutions to TFBM-driven Langevin equations, and \cite{chevillard:2017} on turbulence modeling based on regularized colored noise. In this paper, we approach stochastic differential systems from the dual perspective of \textit{integration}. For the purpose of stochastic analysis, TFLPs are finite variation processes when $d>1/2$ (see Proposition \ref{semimartingaleTFLP}), and hence integration with respect to these processes can be defined pathwise in the usual Stieltjes manner. However, like FLP, when $-1/2<d<1/2$ TFLPs may not be finite variation processes, or even semimartingales (Proposition \ref{p:not_a_semimg} and Remark \ref{r:d=1/2}). For this parameter range, we construct the theory of Wiener-like integrals with respect to these processes. Our approach follows the seminal work \cite{PipirasTaqqu} for FBM, later extended in \cite{MeerschaertsabzikarSPA} to TFBM. Whereas the integration theory with respect to FBM draws upon classical fractional derivatives \cite{FCbook,oldham,Samko}, we put forward a framework for TFLPs based on tempered fractional derivatives~\cite{cartea:del-castillo-negrete:2007,baeumer:meerschaert:2010}. Tempering produces a more tractable mathematical object, and can be made arbitrarily light, so that the resulting operators approximate the fractional derivative to any desired degree of accuracy over compact intervals.

We focus on integration with respect to TFLP \textit{II} (denoted $S^{I\!I}_{d,\lambda}$, $\lambda >0$), since the claims for TFLP are analogous to those for TFBM (see Remark \ref{r:stoch_integ_wrt_TFLP}). Our construction follows from characterizing the natural inner product spaces of integrands ${\mathcal A}_1$ and ${\mathcal A}_2$ (see \eqref{eq:A1 star class} and \eqref{eq:A2class}), which are associated with the memory parameter ranges $-1/2 < d < 0$ and $d > 0$, respectively. In particular, we show that, for TFLP \textit{II}, the phenomenon revealed in \cite{PipirasTaqqu} for FBM resurfaces in the context of tempered fractional L\'{e}vy-type stochastic integration. In other words, for $-1/2 < d < 0$, ${\mathcal A}_1$ and the space of stochastic integrals $\overline{\rm Sp}(S^{I\!I}_{d,\lambda})$ are isometric. As a consequence, every random variable in $\overline{\rm Sp}(S^{I\!I}_{d,\lambda})$ with $-1/2<d<0$ can be written as an integral of a single deterministic function with respect to the stochastic process $S^{I\!I}_{d,\lambda}$ (see Theorems
\ref{thm:stochasticcalculus for TFI and TFD and moving general} and \ref{thm:SLRDisometric2}). However, for $d > 0$, our results show that ${\mathcal A}_2$ is isometric only to a subspace of $\overline{\rm Sp}(S^{I\!I}_{d,\lambda})$ (see Theorems \ref{thm:stochasticcalculus for TFI and moving general} and \ref{thm:SLRDisometric}).

The paper is organized as follows. Section \ref{sec2} contains the definitions and fundamental properties of TFLPs, where Sections \ref{s:TFLP} and \ref{s:TFLPII} pertain to TFLP and TFLP \textit{II}, respectively. In Section \ref{Sec3}, we first show that TFLP and TFLP $I\!I$ are semimartingales for $d>1/2$ and then construct the theory of stochastic integration with respect to these processes for $-1/2<d<1/2$. In Section % \ref{app}, we carry out the data analysis. In Section
 \ref{s:conclusion}, we sum up the conclusions and discuss open problems as well as future research directions. All proofs can be found in the Appendix.  %\BCBcomment{We seem to not discuss the case $d=1/2$ anywhere in the paper, though it appears in some statements implicitly.} \GDcomment{Indeed. Is there a simple way to address the issue? For example, can we add a remark describing what generally happens when $d = 1/2$. Same for $d = 0$, which is probably a less interesting case.} \BCBcomment{I think $d=1/2$ may be a special case and potentially merits special attention, or a remark precluding it (``to avoid technicalities...").}

\section{Moving average representation }\label{sec2}

Recall that a L\'{e}vy process is a stochastically continuous process with stationary and independent increments that starts at zero and has c\`{a}dl\`{a}g sample paths a.s.\ ~\cite{Sato}. Throughout this paper, L\'{e}vy noise plays the role that Brownian noise plays in a Gaussian framework. So, let $L =\{L(t)\}_{t\in\rr}$ be a two-sided L\'{e}vy process constructed by taking two independent copies $L_1=\{L_1(t)\}_{t\geq 0}$ and $L_2=\{L_2(t)\}_{t\geq 0}$ of a L\'{e}vy process and by setting
\begin{equation}\label{e:L(t)}
 L(t):= L_1(t){\bf 1}_{[0,\infty)}(t)-L_2((-t)-){\bf 1}_{(-\infty,0)}(t).
\end{equation}
Hereinafter, we assume $L$ as in \eqref{e:L(t)} satisfies the following condition.

\medskip

\noindent {\bf Condition $L$}: The L\'{e}vy process $L$ in \eqref{e:L(t)} is centered ($\E[ L(1)]=0$) and contains no Brownian component. The distribution of $L$ is uniquely determined by the characteristic function (ch.f.) $\mathbb{E}[\exp{i\theta L(t)}]=\exp\{t\psi(\theta)\}$ for $t\geq 0$, where
\begin{equation}\label{eq:psi defn}
\psi(\theta)=\int_{\rr}(e^{i\theta x}-1-i\theta x)\nu(dx),\quad\quad \theta\in\rr.
\end{equation}
In \eqref{eq:psi defn}, $\nu(dx)$ is called the L\'{e}vy measure of $L$, i.e.,
\begin{equation*}
\nu(\{0\}) =0, \qquad \int_\rr (|x|^2 \wedge 1) \nu(dx)<\infty.
\end{equation*}
Moreover, $\nu(dx)$ is assumed to be such that $\int_{|x|>1}x^{2}\nu(dx)<\infty$, i.e., $\mathbb{E}[(L(t))^2]=t\mathbb{E}[(L(1))^2]=t\int_{\rr}|x|^{2}\nu(dx)<\infty$ for all $t\in\rr$.\\
\medskip

We recall the following classical result for later reference. It provides the conditions for the existence, in the $L^2(\Omega)$ sense, of Wiener-like stochastic integrals with respect to L\'{e}vy noise.
\begin{prop}\label{prop:RajputRosinski}\cite{RajputRosinski,KluppelbergMatsui}
Let $f:\rr\times \rr\to\rr$ be a measurable function. Let $L$ be a L\'{e}vy process such that $\mathbb{E}[L(1)]=0$ and $\mathbb{E}[(L(1))^2] < \infty$. For $t\in\rr$, let $f_t(\cdot)\in L^2(\rr)$. Then, the stochastic integral $S(t):=\int_{\rr}f_{t}(u) dL(u)$ exists in the $L^2(\Omega)$ sense for any $t \in \R$. Furthermore, for $t\in\rr$, $\mathbb{E}[S(t)] = 0$. The isometry
\begin{equation}\label{eq:isometry}
\mathbb{E}[(S(t))^2]=\mathbb{E}[(L(1))^2]\|f_t\|^{2}_{L^{2}(\rr)}, \quad t \in \R,
\end{equation}
also holds, as well as the relation
\begin{equation}
\widetilde{\Gamma}(s,t)={\rm cov}(S(s),S(t))=\mathbb{E}[(L(1))^2]\int_{\rr}f_{s}(u)f_{t}(u)du, \quad s,t \in \R,
\end{equation}
Moreover, the ch.f.\ of $S(t_1),\ldots, S(t_m)$ for $-\infty < t_1< \ldots < t_m < \infty$ is given by
\begin{equation}
\mathbb{E}\Big[\exp\Big\{\sum_{j=1}^{m}i\theta_{j}S(t_j)\Big\}\Big]=\exp\Big\{\int_{\rr}\psi\Big(\sum_{j=1}^{m} \theta_{j}f_{t_j}(s) \Big)ds\Big\}
\end{equation}
for $\theta_j\in\rr, j=1,2,\ldots,m$, where $\psi$ is given by \eqref{eq:psi defn}.
\end{prop}

\subsection{Tempered fractional L\'{e}vy processes of the first kind}\label{s:TFLP}

In this section, we introduce and study tempered fractional L\'{e}vy process of the first kind. We start with its definition.
\begin{defn}\label{defTFLP}
{ \rm Let $L =\{L(t)\}_{t\in\rr}$ be the two-sided L\'{e}vy process \eqref{e:L(t)}. Consider the function $(x)_{+}=xI(x>0)$ and set the convention $0^0=0$. Consider the function $g^{I}_{d,\lambda,t}: \R \rightarrow \R$ given by
\begin{equation*}\label{integrand g I}
g^{I}_{d,\lambda,t}(x):=e^{-\lambda(t-x)_{+}}{(t-x)_{+}^{d}}-e^{-\lambda(-x)_{+}}{(-x)_{+}^{d}}.
\end{equation*}
For any $d>-\frac{1}{2}$ and $\lambda> 0$, the stochastic process
\begin{equation}\label{eq:TFLPdef}
S^{I}_{d,\lambda}(t):=\frac{1}{\Gamma(1+d)} \int_{\rr}g^{I}_{d,\lambda,t}(x) dL(x), \quad t\in\R,
\end{equation}
is called a {\it tempered fractional L\'{e}vy process of the first kind} \textrm{(TFLP)}.}
\end{defn}
The kernel function $g^{I}_{d,\lambda,t}(x)$ is square integrable over $\R$. Hence, by Proposition \ref{prop:RajputRosinski}, the stochastic integral in \eqref{eq:TFLPdef} exists in the $L^2(\Omega)$ sense for any $t \in \R$.

The class of stochastic processes given by Definition \ref{defTFLP} is closely related to a number of other frameworks. When $-\frac{1}{2}<d<\frac{1}{2}$ and tempering is eliminated ($\lambda = 0$), the expression on the right-hand side of \eqref{eq:TFLPdef} is the classical FLP. If $d=0$ (and $\lambda>0$), then $S^{I}_{0,\lambda}(t)$ is called a L\'{e}vy Ornstein-Uhlenbeck (OU) process (\cite{Sato}, Section 3.17). If $dL(x)$ in \eqref{eq:TFLPdef} is replaced with a Gaussian random measure, the resulting process is a TFBM.

Hereinafter, for $S^{I}_{d,\lambda}$ we assume $d\neq0$ (and $\lambda > 0$) unless otherwise stated. Note also that, for any $s,t\in\rr$, the integrand \eqref{integrand g I} satisfies $g^{I}_{d,\lambda,s+t}(s+x)-g^{I}_{d,\lambda,s}(s+x)=g^{I}_{d,\lambda,t}(x)$, and hence one can show that TFLP has stationary increments. In the next proposition, we provide the covariance structure of TFLP.

\begin{prop}\label{prop:covariance_TFLP}
A TFLP $S^{I}_{d,\lambda}$ (see \eqref{eq:TFLPdef}) has the covariance function
\begin{equation}\label{eq:TFLPacvf}
{\rm Cov}\left[S^{I}_{d,\lambda}(t),S^{I}_{d,\lambda}(s)\right]=\frac {\mathbb{E}[(L(1))^2]}{2\Gamma(1+d)^2}\Big\{|t|^{1+2d}C^{2}_{d,\lambda,|t|}
+|s|^{1+2d}C^{2}_{d,\lambda,|s|}
-|t-s|^{1+2d}C^{2}_{d,\lambda,|t-s|}\Big\}
\end{equation}
for any $s,t\in\rr$. %, where $H= d+ \frac{1}{2}$ \GDcomment{We are not using $H$ above}.
In \eqref{eq:TFLPacvf}, %Here
\begin{equation}\label{eq:CtDef}
C^{2}_{d,\lambda,|t|}=\frac{2\Gamma(1+2d)}{(2\lambda |t|)^{1+2d}}-\frac{2\Gamma(1+d)}{\sqrt{\pi }}\Big(\frac{1}{2\lambda |t|}\Big)^{\frac{1}{2}+d} K_{\frac{1}{2}+d}(\lambda |t|),
\end{equation}
for $t\neq 0$, and we define $C_{d,\lambda,0}^{2}=0$. In \eqref{eq:CtDef}, $K_{\nu}(z)$ is the modified Bessel function of the second kind, which is given by
$$
K_{\nu}(z) = \int^{\infty}_0 e^{-z (\frac{e^{-t}+e^{t}}{2})} \frac{e^{-\nu t}+e^{\nu t}}{2} dt, \quad z > 0, \quad \nu \in \R.
$$
Moreover,
\begin{equation}\label{e:limitbehaviorofTFLP}
\lim_{t\to\infty} {\rm Var} \big[ S^{I}_{d,\lambda} (t) \big] = \frac{2 \E( L(1)^2 ) \Gamma(1+2d)}{ \Gamma(1+d)^2 (2\lambda)^{1+2d} }.
\end{equation}
\end{prop}
It is well known that the variance of FLP is divergent \cite{Marquardt}. Remarkably, expression \eqref{e:limitbehaviorofTFLP} shows that the variance of TFLP stays finite in the large scale limit (cf.\ \cite{boniece:didier:sabzikar:2018}, Proposition A.1).

\begin{rem}\label{TFBMBSSdifference1}
{\rm
Let $L =\{L(t)\}_{t\in\rr}$ be the two-sided L\'{e}vy process \eqref{e:L(t)}. Then, a L\'{e}vy semistationary process (LSS; see \cite{BarndorffNielsen1,barndorff-nielsen:2016}) is defined by the stochastic integral representation
\begin{equation}
Y(t)= \mu+ \int_{-\infty}^{t} g(t-s)\sigma(s) dL(s) + \int_{-\infty}^{t} q(t-s)a(s) ds,
\end{equation}
where $\sigma$ and $a$ are stochastic processes, and $g$ and $q$ are deterministic kernels with $g(t) = h (t) = 0$ for $t \leq 0$. Although LSS instances associated with gamma kernels ($g(x)= x^{d-1}e^{-\lambda x}$) and TFLP both display a tempering component, the two processes are generally quite different. In particular, the former may be stationary, while the latter is always nonstationary.
}
\end{rem}

In the next proposition, we establish a stochastic integral representation of TFLP as an improper Riemann integral for the parameter range $d > 0$. The result is then used in part $(a)$ of the subsequent theorem to construct a H\"older--continuous modification of TFLP.
%%%%%%%%%%%%%%%%%%%
%The next Theorem shows TFLP has a continuous modification by showing that the stochastic integral in
%\eqref{eq:TFLPdef} can be written as an improper Riemann integral.\GDcomment{I don't understand the claim. The improper Riemann integral is used in the proof of the proposition. Moreover, continuity only holds for $0 < d < 1/2$.}
%%%%%%%%%%%%%%

\begin{prop}\label{p:modification_Theorem_TFLP1} Let $S^{I}_{d,\lambda}= \{ S^{I}_{d,\lambda}(t) \}_{t\in\rr}$ be a TFLP (see \eqref{eq:TFLPdef}) with $d > 0$.  Then, for all $t\in\rr$, there exists a modification of $S^{I}_{d,\lambda}(t)$ which is equal to the improper Riemann integral
 \begin{equation}\label{e:SI_improper_Riemann}
\begin{split}
S^{I}_{d,\lambda}(t) &=  \frac{1}{\Gamma(d)}
 \int_{\rr}\Big( e^{-\lambda(t-x)_{+}}{(t-x)_{+}^{d-1}}-e^{-\lambda(-x)_{+}}{(-x)_{+}^{d-1}}\Big)\ L(x)\ dx\\
&\qquad\qquad - \frac{\lambda}{\Gamma(d+1)} \int_{\rr}\Big( e^{-\lambda(t-x)_{+}} (t-x)_{+}^{d} - e^{-\lambda(-x)_{+}} (-x)_{+}^{d} \Big) \ L(x)\ dx.
\end{split}
\end{equation}
In particular, the process \eqref{e:SI_improper_Riemann} is continuous in $t$.
\end{prop}

The following theorem is our main result on the sample path properties of TFLP. Note that the statement in $(a)$ is slightly stronger than the one usually obtained in the framework of the Kolmogorov-$\breve{\textnormal{C}}$entsov criterion. %\GDcomment{What prevents us from assuming $d > 1/2$? Same for TFLP \textit{II}.}
\begin{thm}\label{t:TFLP_sample_path}
Let $S^{I}_{d,\lambda} = \{ S^{I}_{d,\lambda}(t) \}_{t\in\rr}$ be a TFLP (see \eqref{eq:TFLPdef}).
\begin{itemize}
\item [(a)] If $0<d\leq\frac{1}{2}$, then %for every $0 < \gamma < d$
there exists a locally $d$-H\"{o}lder continuous modification of $S^{I}_{d,\lambda}$. That is, for $T > 0$,
\begin{equation}\label{TFLPsamplepath}
\mathbb{P} \Bigg[ \omega: \sup_{ 0<|s-t|<k_{T}(\omega), |s|\leq T, |t|\leq T  } \Big( \frac{| S^{I}_{d,\lambda}(t) - S^{I}_{d,\lambda}(s) | }{|s-t|^d}  \Big)\leq C       \Bigg] =1,
 \end{equation}
where $k_{T}(\omega)$ is an almost surely positive random variable and $C >0$.
%\noi(b)
\item [(b)] If $-\frac{1}{2} < d <0$ and $L$ has symmetric finite-dimensional distributions, then $S^{I}_{d,\lambda}$ has discontinuous and unbounded sample paths with positive probability.
\end{itemize}
\end{thm}

Next, we turn to the increment process of TFLP. Starting from a TFLP $S^{I}_{d,\lambda}$, the stationary process \textit{tempered fractional L\'{e}vy noise of the first kind} (TFLN) is naturally defined as
\begin{equation}\label{eq:TFLNdef}
X^{I}_{d,\lambda}(t) := S^{I}_{d,\lambda}(t+1)-S^{I}_{d,\lambda}(t), \quad t\in\rr.
\end{equation}
It follows readily from \eqref{eq:TFLPdef} that TFLN has the moving average representation
\begin{equation}\label{eq:TFGNmoving}
X^{I}_{d,\lambda}(t)=\frac{1}{\Gamma(d+1)}
\int_{\rr} \big[ e^{-\lambda(t+1-x)_{+}} (t+1-x)_{+}^{d} - e^{-\lambda(t-x)_{+}} (t-x)_{+}^{d} \big]\ dL(x).
\end{equation}
In the following proposition, we characterize the behavior of the covariance of TFLN over large lags. In particular, TFLN is semi-LRD in the sense of \eqref{e:semi-LRD} with $\delta = d > -1/2$.
\begin{prop}\label{prop:TFLN}
Let $X^{I}_{d,\lambda}=\{X^{I}_{d,\lambda}(t)\}_{t\in\R}$ be a TFLN (see \eqref{eq:TFGNmoving}). Let $\gamma^{I}(h)=\E [X^{I}_{d,\lambda}(0) X^{I}_{d,\lambda}(h)] $ be its covariance function and let $h^{I}(\omega)$ be its spectral density. Then,
\begin{itemize}
%\noi(a)
\item [(a)] as $h\to\infty$,
\begin{equation}\label{e:gammaI(h)_decay_semi-LRD}
\gamma^{I}(h) \sim C e^{-\lambda h} h^{d} ,
\end{equation}
%where $d > -\frac{1}{2}$ and
where $C = C(d,\lambda)= - \frac{\E [L(1)^2] \lambda^2}{ \Gamma(d+1) (2\lambda)^{d+1}} $; %is a constant depending on the parameters $d$ and $\lambda$;
%\noi(b) $\{ X(t) \}_{t\in\mathbb{R}}$ has negative memory.

%\noi(b)
\item [(b)] for $\omega\in \rr$,
\begin{equation}\label{e:TFLN_density}
h^{I}(\omega) = \frac{1}{2\pi} \frac{(1-\cos{\omega} )}{ (\lambda^2 + \omega^2)^{d+1} }.
\end{equation}
\end{itemize}
\end{prop}

\subsection{Tempered fractional L\'{e}vy processes of the second kind}\label{s:TFLPII}

In this section, we introduce and study tempered fractional L\'{e}vy process of the second kind. We start with its definition.
\begin{defn}\label{e:defTFLPII}
{\rm Let $L =\{L(t)\}_{t\in\rr}$ be the two-sided L\'{e}vy process \eqref{e:L(t)} and consider the function $g^{I\! I}_{d,\lambda,t}: \R \to \R $ given by
\begin{equation*}\label{hdef0}
g^{I\! I}_{d,\lambda,t}(y)=(t-y)_+^{d} \ e^{-\lambda (t-y)_+} - (-y)_+^{d} \ e^{-\lambda (-y)_+}
+ \lambda \int_{0}^{t} (s-y)_{+}^{d} \ e^{-\lambda(s-y)_{+}} \ ds.
\end{equation*}
For any $d> -\frac 12$ and $\lambda> 0$, the stochastic process
\begin{equation}\label{eq:defTFLP second}
S^{I\!I}_{d,\lambda}(t):=\frac{1}{\Gamma(d+1)} \int_{\rr} g^{I\! I}_{d,\lambda,t}(y)\ dL(y),\quad  t \in \R,
\end{equation}
is called a {\it tempered fractional L\'{e}vy process of the second kind} (TFLP \textit{II})}.
\end{defn}
By Proposition \ref{prop:RajputRosinski}, $S^{I\!I}_{d,\lambda}(t)$ is well defined in the $L^2(\Omega)$ sense for any $t \in \R$, since $g^{I\! I}_{d,\lambda,t}(y)$ is square integrable (see Lemma \ref{lem:g squar integrable}).

As with \eqref{eq:TFLPdef}, the class of stochastic processes given by \eqref{eq:defTFLP second} is closely related to other frameworks. When $-\frac{1}{2}<d<\frac{1}{2}$ and tempering is eliminated ($\lambda=0$), the process $S^{I\!I}_{d,0}(t)$ also reduces to FLP. If $dL(x)$ in \eqref{eq:defTFLP second} is replaced with a Gaussian random measure, the resulting process is a TFBM $I\!I$.

Hereinafter, for $S^{I\!I}_{d,\lambda}$ we assume $d\neq0$ (and $\lambda > 0$), unless otherwise stated.

In the following proposition, we express the covariance function $\E[S^{I\!I}_{d,\lambda}(t)S^{I\!I}_{d,\lambda}(s)]$ of TFLP $I\!I$ when $d > 0$.
\begin{prop}\label{prop:THPcovariance}
For $d>0$, a TFLP $I\!I$ (see \eqref{eq:defTFLP second}) has covariance function
\begin{equation}\label{eq:THPacvf}
{\rm Cov}\left[S^{I\! I}_{d,\lambda}(t),S^{I\! I}_{d,\lambda}(s)\right]=\frac{\E[L(1)^2]}{\sqrt{\pi}\Gamma(d)(2\lambda)^{d-\frac{1}{2}}}\int_{0}^{t}\int_{0}^{s}
|u-v|^{d-\frac{1}{2}}K_{d-\frac{1}{2}}(\lambda|u-v|)dv\ du
\end{equation}
for any $s,t\in\rr$.
\end{prop}

\begin{rem}
{\rm In the parameter range $-1/2<d<0$, the covariance of TFLP $I\!I$ can be found by first developing $\E[(S^{I\!I}_{d,\lambda}(t))^2]$ via the isometry property \eqref{eq:isometry}, and then applying the elementary formula $ab=\frac{1}{2}(a^2+b^2-(a-b)^2)$ as well as the stationary increments property.  However, the final formula for $\E[(S^{I\!I}_{d,\lambda}(t))^2]$ involves several integral expressions, and consequently so does the formula for $\E[S^{I\!I}_{d,\lambda}(t)S^{I\!I}_{d,\lambda}(s)]$.  For brevity and clarity of exposition, we opt for not including it here.}\end{rem}

The next proposition is the analog for TFLP $I\!I$ of Proposition \ref{p:modification_Theorem_TFLP1}. It shows that TFLP $I\!I$ has a modification that can be written as an improper Riemann integral.  %\BCBcomment{We don't use the Riemann integral to get the holder regularity for this one.}
\begin{prop}\label{modification_Theorem_TFLP2}
Let $S^{I\!I}_{d,\lambda} = \{S^{I\!I}_{d,\lambda}(t)\}_{t \in \R}$ be a TFLP $I\!I$ (see \eqref{eq:defTFLP second}) with $d>0$. Then, for all $t\in\rr$, there exists a modification of $S^{I\!I}_{d,\lambda}(t)$ which is equal to the improper Riemann integral
\begin{equation}\label{e:SII_improper_Riemann}
\begin{split}
S^{I\!I}_{d,\lambda}(t) &= \frac{1}{d\Gamma(d-1)} \int_{\rr}\int_0^t  e^{-\lambda(s-x)_{+}} (s-x)_{+}^{d-2} ds L(x) dx\\
&- \frac{\lambda}{\Gamma(d-1)} \int_{\rr}\int_0^t  e^{-\lambda(s-x)_{+}} (s-x)_{+}^{d-1} ds L(x) dx.
\end{split}
\end{equation}
In particular, the process \eqref{e:SII_improper_Riemann} is continuous in $t$.
\end{prop}

The following theorem is our main result on the sample path properties of TFLP $I\!I$.
\begin{thm}\label{t:TFLPI bounds}Let  $S^{I\!I}_{d,\lambda}= \{S^{I\!I}_{d,\lambda}(t)\}_{t \in \R}$ be a TFLP $I\!I$ (see \eqref{eq:defTFLP second}).
\begin{itemize}
\item [(a)] If $0<d\leq\frac{1}{2}$, then for every $0 < \gamma < d$, there exists a locally $\gamma$-H\"{o}lder continuous modification of $S^{I\! I}_{d,\lambda}$. That is, for $T > 0$,
\begin{equation}\label{TFLPsamplepath}
\mathbb{P} \Bigg[ \omega: \sup_{ 0<|s-t|<k_{T}(\omega), |s|\leq T, |t|\leq T  } \Big( \frac{ |S^{I\! I}_{d,\lambda}(t) - S^{I\! I}_{d,\lambda}(s)  |}{|s-t|^\gamma}\Big) \leq C       \Bigg] =1,
 \end{equation}
where $k_{T}(\omega)$ is an almost surely positive random variable and $C >0$.
\item [(b)] If $-\frac{1}{2} < d <0$ and $L$ has symmetric finite-dimensional distributions, %\GDcomment{To be confirmed that this is what is meant by ``symmetric".},
 then $S^{I\! I}_{d,\lambda}$ has discontinuous and unbounded sample paths with positive probability.
\end{itemize}
%\noi(b)
\end{thm}

Next, we turn to the increment process of TFLP $I\!I$. Starting from a TFLP $I\!I$ $S^{I\!I}_{d,\lambda}$, the stationary process \textit{tempered fractional L\'{e}vy noise of the second kind} (TFLN $I\!I$) is naturally defined as
\begin{equation}\label{eq:defTFLN second}
X^{I\!I}_{d,\lambda}(t) = S^{I\!I}_{d,\lambda}(t+1)-S^{I\!I}_{d,\lambda}(t), \quad t \in \R.
\end{equation}
It follows from \eqref{eq:defTFLP second} that TFLN $I\!I$ has moving average representation
\begin{equation}\label{eq:TFGNmovingII}
X^{I\!I}_{d,\lambda}(t)=\frac{1}{\Gamma(d)}
\int_{\rr} \int_{t}^{t+1}  (s-y)_{+}^{d} e^{-\lambda(s-y)_{+}} \ ds\ dL(y) .
\end{equation}

The following proposition describes the behavior of the covariance structure of TFLN $I\!I$ over large lags. In particular, the proposition shows that TFLN$I\!I$ is semi-LRD in the sense of \eqref{e:semi-LRD} with $\delta = d - 1 > -3/2$. In the Fourier domain, it shows that the spectral density is of the Von K\'{a}rm\'{a}n type (cf.\ Figure \ref{fig:vonkarman}). In the statement of the proposition, we make use of the following notation: given two real-valued functions $f(t)$, $g(t)$ on $\mathbb{R}$, we write $f(t)\asymp g(t)$ if $C_1\leq |{f(t)}/{g(t)}|\leq C_2$ for all $t>0$ sufficiently large, for some $0<C_1<C_2<\infty$.
\begin{prop}\label{p:TFLNII}
Let $X^{I\!I}_{d,\lambda} = \{X^{I\!I}_{d,\lambda}(t)\}_{t \in \R}$ be a TFLN $I\!I$ (see \eqref{eq:defTFLN second}). Let $\gamma^{I\!I}(h)=\E [X^{I\!I}_{d,\lambda}(0) X^{I\!I}_{d,\lambda}(h)] $, $h \in \R$, be its covariance function, and let $\{h^{I\!I}(\omega)\}_{\omega \in\R}$ be its spectral density. Then,
%\noi(a)
\begin{itemize}
\item [(a)] as $h\to\infty$,
\begin{equation}\label{e:gamma(h)_semi-LRD_TFLNII}
\gamma^{I\!I}(h) \asymp e^{-\lambda h } h^{ d -1 };
\end{equation}%\item [(b)] $\{ X(t) \}_{t\in\mathbb{R}}$ has short memory.
%\noi(b)
\item [(b)] for $\omega\in \rr$,
\begin{equation*}
h^{I\!I}(\omega) = \frac{1}{2\pi} \frac{(1-\cos{\omega} )}{ \omega^2\ (\lambda^2 + \omega^2)^{d} }. %\frac{\sigma^2}{2\pi}
\end{equation*}
\end{itemize}
\end{prop}

As a preparation for the next section -- on stochastic integration --, we conclude this section by constructing subclasses of TFLPs that are not semimartingales. Note that, in all cases, the memory parameter is taken in the range $d \in (-1/2,1/2)$.

So, let $(\mathcal{F}^{L,\infty}_{t})_{t\geq 0}$ be the smallest filtration such that $\sigma\Big( L(s): -\infty < s \leq t  \Big)\subseteq (\mathcal{F}^{L,\infty}_{t})_{t\geq 0}$ for all $t\geq0$.% Also, let $(\mathcal{F}^{L}_{t})_{t\geq 0}$ be the natural filtration of $\{ L(t) \}_{t\geq 0}$.%The proposition below is analogous to Theorem 4.5 in \cite{Marquardt}.
\begin{prop}\label{p:not_a_semimg}
Let $1<\alpha<2$, and choose $d\in(-1/2,1/2)\setminus\{0\}$ such that $d+\frac{1}{\alpha}\in (0,1)$. For such $d$, let $S^{*}_{d,\lambda} = \{ S^{*}_{d,\lambda}(t) \}_{t\in\rr}$ be a either a TFLP  \eqref{eq:TFLPdef} or TFLP \textit{II}  \eqref{eq:defTFLP second} with $\nu(dx)=h(x)dx$, where
\begin{equation}\label{h(x)_xto0}
 h(x)\sim |x|^{-1-\alpha}, \quad x\to0.
\end{equation}
Then, $S^{*}_{d,\lambda}$ is not a $(\mathcal{F}^{L,\infty}_{t})_{t\geq 0}$-semimartingale.
\end{prop}

\begin{figure}[h!]
    \centering
    \begin{minipage}{.5\linewidth}
        \centering
        \includegraphics[width=\linewidth]{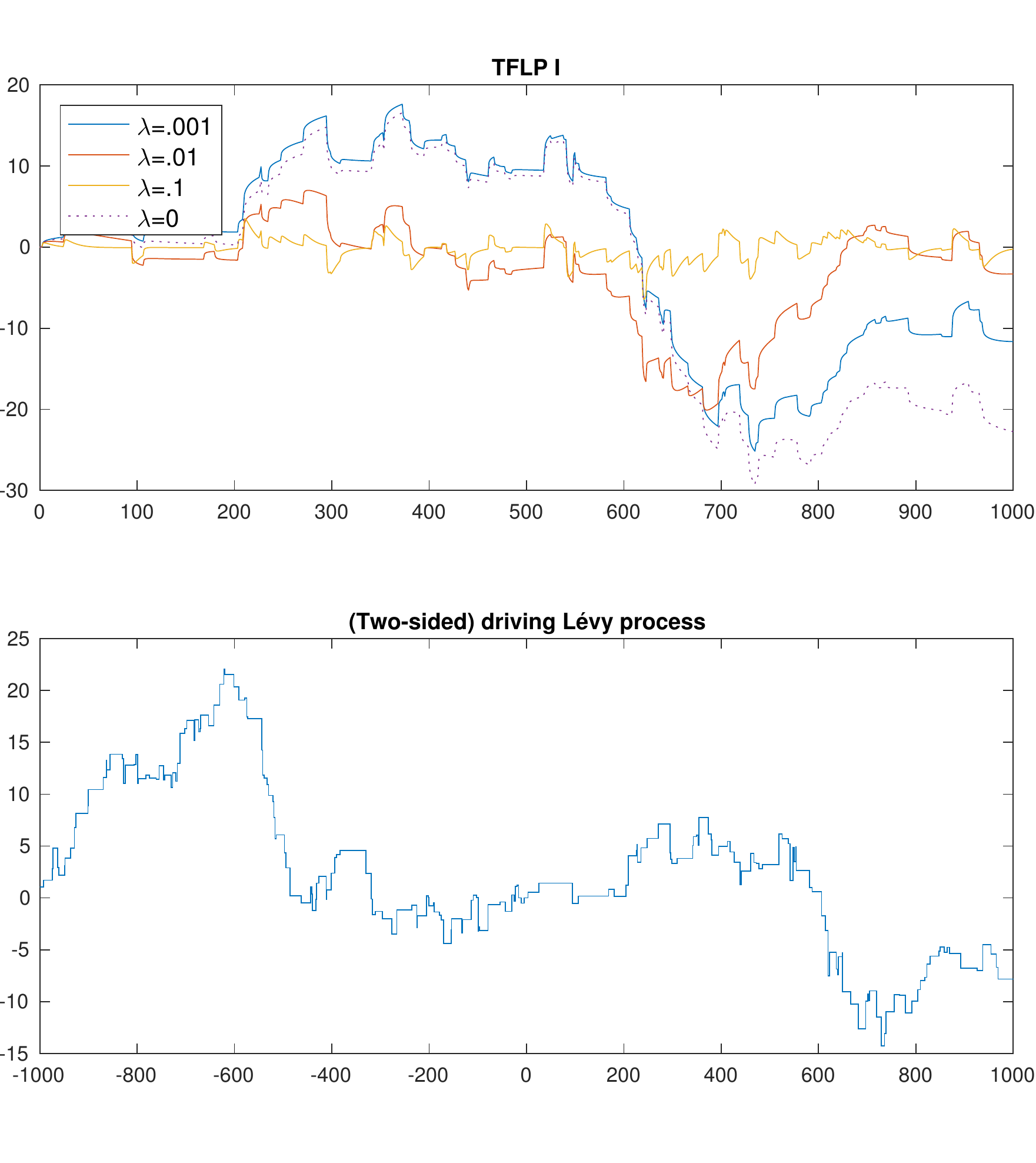}
        \caption*{Above: compound Poisson driving noise}
        \label{fig:TFLP_paths1}
    \end{minipage}%
    \begin{minipage}{.5\linewidth}
        \centering
        \includegraphics[width=\linewidth]{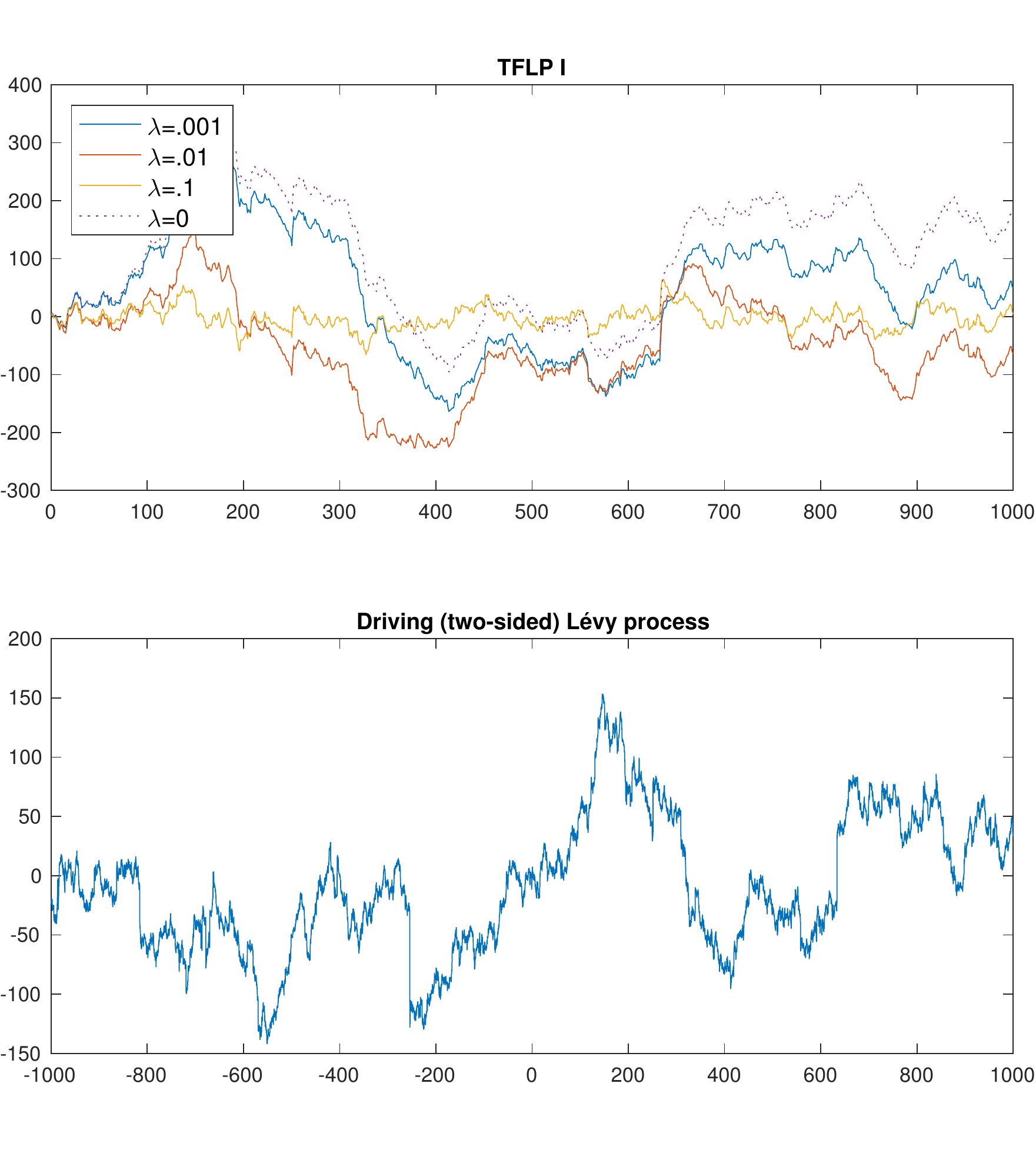}
        \caption*{Above: tempered stable driving noise}
        \label{fig:TFLP_paths2}
    \end{minipage} \caption{\label{fig:TFLPI_paths} \textbf{Simulated paths of $S_{d,\lambda}^I$}. In the top row of figures above, paths with memory parameter $d=1/6$ were generated for $\lambda\in\{0,0.001,0.01,0.1\}$ based on the same corresponding driving process (bottom plots).  The plot on the left uses a compound Poisson driving process with intensity $1$ and uniform $[-1,1]$ jumps.  The figure on the right uses symmetric tempered $\alpha$-stable driving noise (see, e.g., \cite{baeumer:meerschaert:2010} for details on the simulation of such processes) with tempering parameter $\lambda_{\text{noise}}=.01$ and $\alpha=1.65$.  The discontinuous sample paths of the driving processes are displayed as continuous lines for visual clarity (see Example \ref{ex:simulation} on simulation details). The convergence to stationarity effect caused by tempering is more visible for larger values of $\lambda$.}
%       }
\end{figure}

%\begin{center}
%\includegraphics[width=.5\linewidth]{sample_paths_alph_16_lam_005.eps}
%\end{center}
\begin{figure}[h!]
    \centering
    \begin{minipage}{.5\linewidth}
        \centering
        \includegraphics[width=\linewidth]{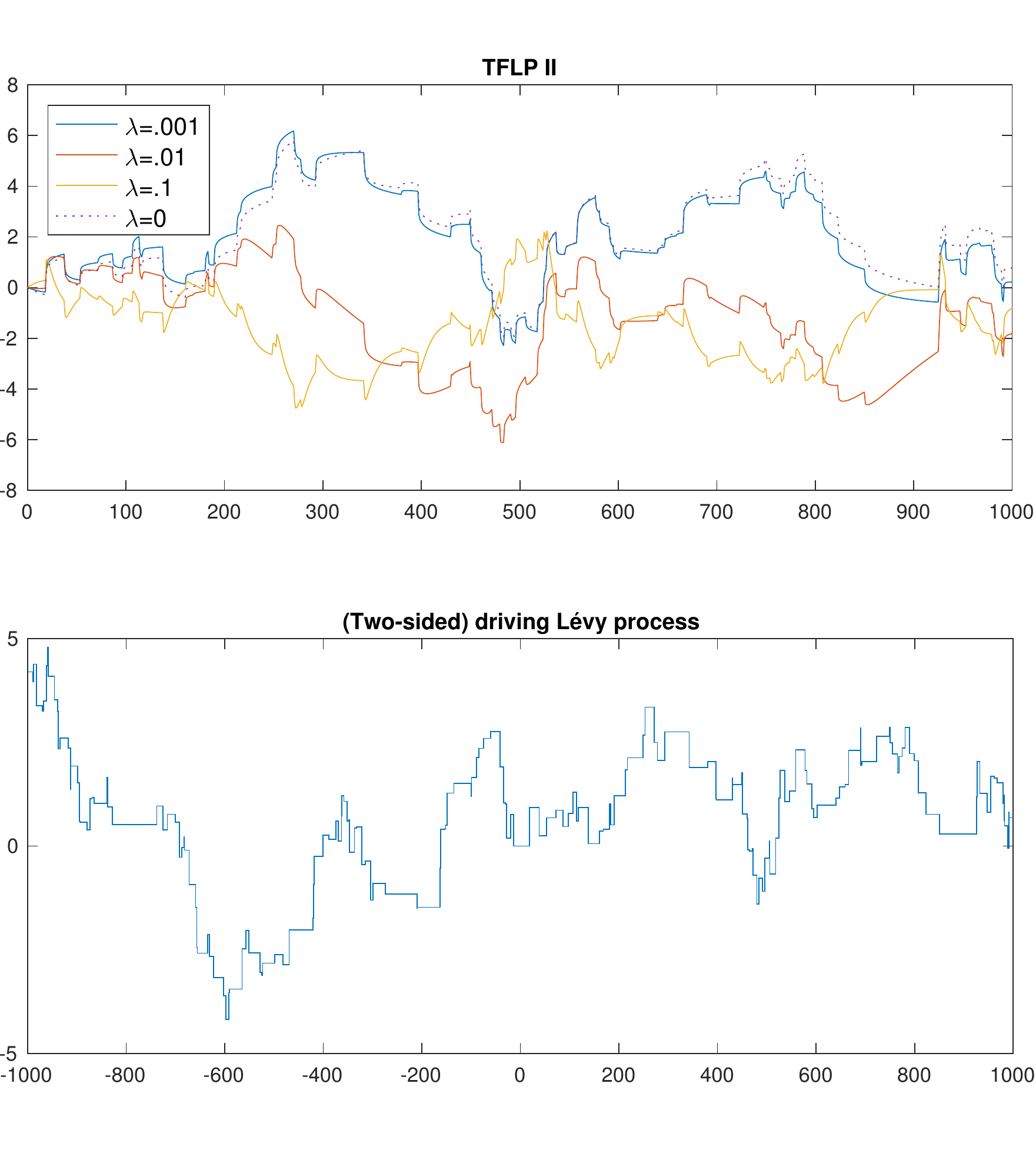}
        \caption*{Above: compound Poisson driving noise}
        \label{fig:TFLP_paths1}
    \end{minipage}%
    \begin{minipage}{.5\linewidth}
        \centering
        \includegraphics[width=\linewidth]{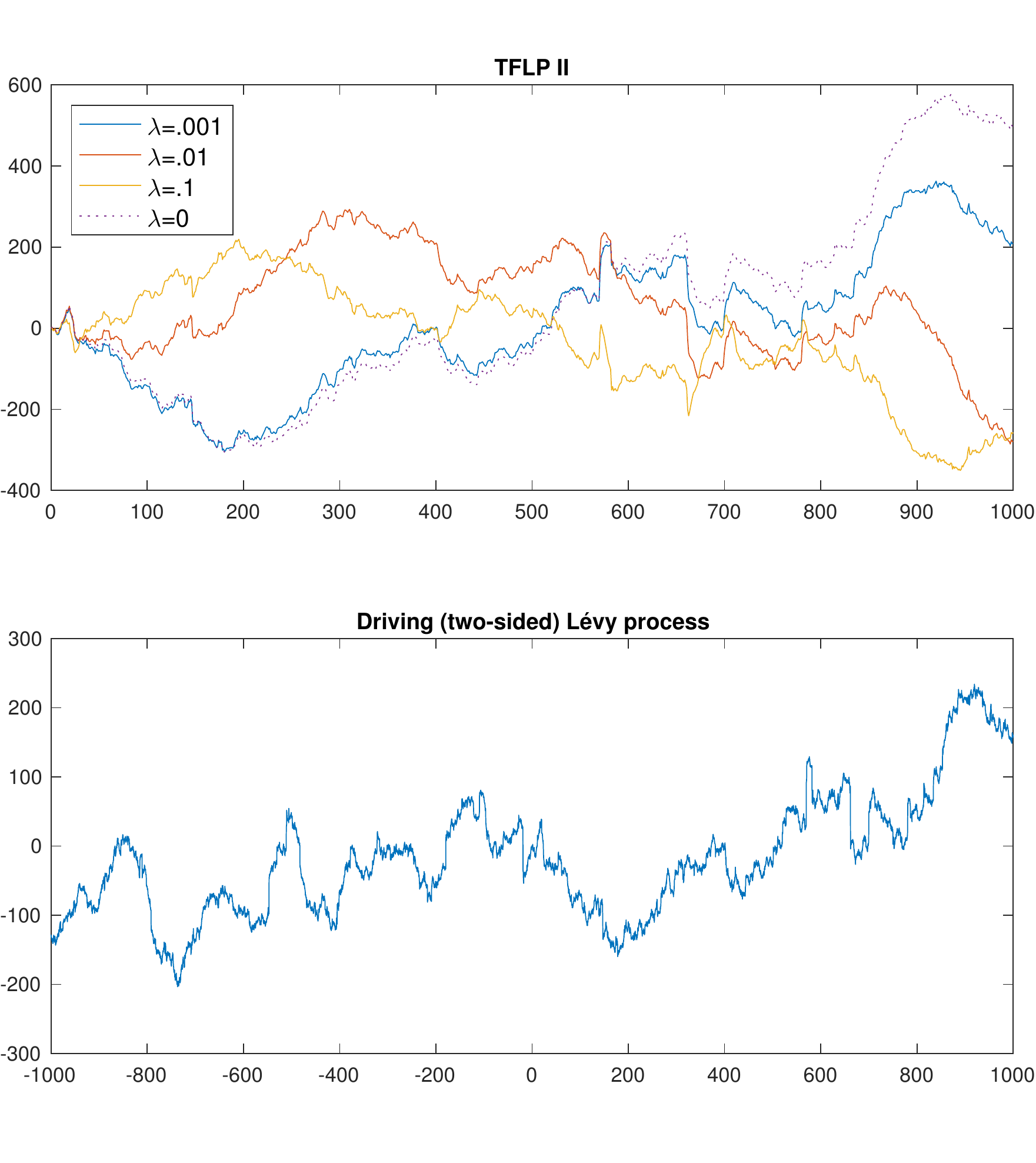}
        \caption*{Above: tempered stable driving noise}
        \label{fig:TFLP_paths2}
    \end{minipage} \caption{\label{fig:TFLPII_paths} \textbf{Simulated paths of $S_{d,\lambda}^{I\!I}$}. In the top row of figures above, paths with memory parameter $d=1/6$ were generated for $\lambda\in\{0,0.001,0.01,0.1\}$ based on the same corresponding driving process (bottom plots). The plot on the left uses a compound Poisson driving process with intensity $1$ and uniform $[-1,1]$ jumps.  The figure on the right uses symmetric tempered $\alpha$-stable driving noise (see, e.g., \cite{baeumer:meerschaert:2010} for details on the simulation of such processes) with tempering parameter $\lambda_{\text{noise}}=.01$ and $\alpha=1.65$. The discontinuous sample paths of the driving processes are displayed as continuous lines for visual clarity (see Example \ref{ex:simulation} on simulation details).   }
%       }
\end{figure}

%\begin{rem}
%{\rm Analogous to \cite[Theorem 4.1]{BenassiCohenIstas1} and  \cite[Theorem 4.5]{Marquardt}, if $S^*_{d,\lambda}$ is taken with L\'evy measure and $\alpha,d$ as above, under the condition $h(x)\leq C |x|^{-1-\alpha}$ for all $x\in \R$, it can be shown that $S^*_{d,\lambda}$ displays local self-similarity with Hurst parameter $H=d+1/\alpha$. In particular, $S^*_{d,\lambda}$ converges to a linear fractional stable motion with parameters $\alpha,H$, i.e., for every fixed $t$,  as $\epsilon\to 0$,
%$$
% \Big\{\frac{S^{*}_{d,\lambda}(t+\epsilon x) - S^{*}_{d,\lambda}(t)}{\epsilon^H}\Big\}_{x\in\R}\stackrel{d}{\to} \{Y_H(x)\}_{x\in\R},
%$$
%where $Y_H(x)= \Gamma(d)^{-1}\int_\R \{(x-s)_+^{H-\frac{1}{\alpha}}-(-s)_+^{H-\frac{1}{\alpha}} \} M(ds)$, and $M$ is a symmetric $\alpha$-stable motion.
%}\end{rem} \vspace{-4ex}
%\BCB{\begin{rem}{\rm Among the class of TFLP and TFLP $I\!I$ with L\'evy measure satisfying \eqref{h(x)_xto0} are those driven by tempered stable processes (see, e.g., \cite{kuchler:tappe:2013}), which can be easily simulated (see, e.g., \cite{baeumer:meerschaert:2010}).}
%\end{rem}}
%\BCBGD{\begin{rem}{ When $\nu(\R)<\infty$, and $d\geq 0$, $S^{I}_{d,\lambda}$ $S^{I\! I}_{d,\lambda}$ are always semimartingales, since they are finite variation processes.  This can be established using the same argument as in \cite{Marquardt}, Theorem 4.6.} % \BCBcomment{I think I can characterize precisely when TFLP I/II are semimartingales (or not), along the lines of \cite{Bender}.  In the interest of time I've left this out for now.} \GDcomment{Fine.}}
%\end{rem}}

\begin{example}\label{ex:simulation}
{\rm Figures \ref{fig:TFLPI_paths} and \ref{fig:TFLPII_paths} display simulated sample paths of TFLP and TFLP $I\!I$. The simulation was carried out based on Riemann-Stieltjes sums, in the fashion of \cite[p.\ 89]{MarquardtphD}. Multiple values of the tempering parameter $\lambda$ and two different types of driving L\'{e}vy noise were used as to illustrate the effect of tempering and of distinct non-Gaussian distributions, respectively.}
\end{example}

\begin{rem}\label{r:d=1/2}
{\rm
The argument for showing Proposition \ref{p:not_a_semimg} requires $d\in(-1/2,1/2)\setminus\{0\}$. Whether or not the boundary value $d=1/2$ always gives a semimartingale remains an open question (cf.\ Proposition \ref{semimartingaleTFLP}).}
\end{rem}

\section{Stochastic integration with respect to TFLP and TFLP $I\!I$}\label{Sec3}

In this section, we develop the theory of stochastic integration with respect to TFLPs. Recall that TFLP and TFLP $I\!I$ are both well defined for $d>-1/2$ and $\lambda>0$.

Stochastic integration theory for FBM and FLP is complicated by the fact that they are not semimartingales \cite{PipirasTaqqu,MarquardtphD}. In contrast, as shown in the following proposition, the representations of TFLP and TFLP $I\!I$ as Riemann-Stieltjes integrals imply that they are finite variation processes when $d>1/2$. Consequently, in this parameter range, we can conveniently define integrals
$$
I(f):=\int f(x)dS^*_{d,\lambda}(x), \quad dS^*_{d,\lambda}(x) = dS^I_{d,\lambda}(x) \hspace{2mm} \textnormal{ or } \hspace{2mm}dS^*_{d,\lambda}(x) = dS^{I\!I}_{d,\lambda}(x),
$$
$\omega$-by-$\omega$ as ordinary Stieltjes integrals (see \cite[p.\ 283]{Kallenberg} or \cite[pp.\ 149--150]{karatzas:shreve:2000}).
\begin{prop}\label{semimartingaleTFLP}%\label{semimartingaleTFLPI}
Suppose $d>1/2$, and let $L = \{ L(t)\}_{t\in\rr}$ be the two-sided L\'{e}vy process \eqref{e:L(t)}.
\begin{itemize}
\item [$(i)$] Let $S^{I}_{d,\lambda}=\{S^{I}_{d,\lambda}(t)\}_{t\in \R}$  be a TFLP (see \eqref{eq:TFLPdef}). Then, the process
\begin{equation}\label{e:TFLP_d_g_half}
\left\{\frac{1}{\Gamma(d+1)}\int_{0}^t \int_{-\infty}^s d(s-x)^{d-1}e^{-\lambda(s-x)} -\lambda (s-x)^{d}e^{-\lambda(s-x)}dL(x) ds\right\}_{t\in\R}
\end{equation}
is a version of $S^{I}_{d,\lambda}$. In particular, for such $d$, $S^{I}_{d,\lambda}$ has a.s.\ absolutely continuous paths and hence is a finite variation process.
\item [(ii)] Let $S^{I\!I}_{d,\lambda}= \{S^{I\!I}_{d,\lambda}(t)\}_{t\in \R}$  be a TFLP $I\!I$ (see \eqref{eq:defTFLP second}). Then, the process
\begin{equation}\label{e:TFLPII_d_g_half}
\left\{\frac{1}{\Gamma(d+1)}\int_{0}^t \int_{-\infty}^s d(s-x)^{d-1}e^{-\lambda(s-x)}dL(x) ds\right\}_{t\in\R}
\end{equation}
is a version of $S^{I\!I}_{d,\lambda}$. In particular, for such $d$, $S^{I\!I}_{d,\lambda}$ has a.s.\ absolutely continuous paths and hence is a finite variation process.
\end{itemize}
\end{prop}

Next, we tackle the case
\begin{equation}\label{e:-1/2<d<1/2_lambda>0}
-1/2 <d < 1/2.
\end{equation}
Even though \eqref{e:-1/2<d<1/2_lambda>0} is our focus, whenever applicable we use the larger range interval $d > 0$ instead of $1/2 > d > 0$.

First, we show the connection between tempered fractional processes and tempered fractional calculus. We refer the reader to the appendix for more details on the latter.

\begin{defn}{\rm \label{defn:Tempered fractional integral}
For any $f\in{L}^{p}({\mathbb{R}})$, $1\leq p<\infty$, the positive and negative \textit{tempered fractional integrals} of a function $f:\R \rightarrow \R$ are defined by
\begin{equation}\label{eq:positivetempered fractional integral}
{\mathbb I}^{\kappa,\lambda}_{+} f(y)=\frac{1}{\Gamma(\kappa)}\int_{-\infty}^{+\infty} f(s)(y-s)_{+}^{\kappa-1}e^{-\lambda(y-s)_{+}}ds
\end{equation}
and
\begin{equation}\label{eq:negativetempered fractional integral}
{\mathbb I}^{\kappa,\lambda}_{-} f(y)=\frac{1}{\Gamma(\kappa)}\int_{-\infty}^\infty f(s)(s-y)_{+}^{\kappa-1}e^{-\lambda(s-y)_{+}}ds
\end{equation}
respectively, for any $\kappa>0$ (and $\lambda >0$).
}\end{defn}

Note that, when $\lambda=0$, these definitions reduce to the (positive and negative) Riemann-Liouville fractional integral, which extends the usual operation of iterated integration to a fractional order \cite{FCbook,oldham,Samko}.  When $\lambda=1$, the operator \eqref{eq:positivetempered fractional integral} is called the Bessel fractional integral \cite[Section 18.4]{Samko}.

The inverse operator of the tempered fractional integral is called tempered fractional derivative. For our purposes, we only require derivatives of order $0<\kappa<1$, which simplifies the presentation.

%\BCBcommentnew{Should there be any restriction on the class of functions here?}
\begin{defn}\label{TFDdef}
{\rm
The positive and negative \textit{tempered fractional derivatives} of a function $f:\R\to\R$ are defined as
\begin{equation}\label{eq:temperedfractionalderivativepositive}
{\mathbb{D}}^{\kappa,\lambda}_{+}f(y)={\lambda}^{\kappa}f(y)+\frac{\kappa}{\Gamma(1-\kappa)}\int_{-\infty}^{y}\frac{f(y)-f(s)}{(y-s)^{\kappa+1}}
\,e^{-\lambda(y-s)}\ ds
\end{equation}
and
\begin{equation}\label{eq:temperedfractionalderivativenegative}
{\mathbb{D}}^{\kappa,\lambda}_{-}f(y)={\lambda}^{\kappa}f(y)+\frac{\kappa}{\Gamma(1-\kappa)}\int_{y}^{+\infty}\frac{f(y)-f(s)}
{(s-y)^{\kappa+1}}\,e^{-\lambda(s-y)}\ ds,
\end{equation}
respectively, for any $0<\kappa<1$ (and $\lambda>0$).
}\end{defn}

Note that expressions \eqref{eq:temperedfractionalderivativepositive} and \eqref{eq:temperedfractionalderivativenegative} reduce to the positive and negative Marchaud fractional derivatives if $\lambda=0$ (cf.\ \cite[Section 5.4]{Samko}).

As pointed out in \cite[p.\ 2367]{MeerschaertsabzikarSPA}, tempered fractional derivatives cannot be defined pointwise for all functions $f \in L^p(\R)$. However, $\mathbb{D}^{\kappa,\lambda}_{\pm}f$ is well defined when $f,f' \in L^1(\R)$. For such $f$, the Fourier transform $\mathcal{F}[\mathbb{D}^{\kappa,\lambda}_{\pm}f]$ satisfies  $\mathcal{F}[\mathbb{D}^{\kappa,\lambda}_{\pm}f](\omega)= (\lambda\pm i\omega)^\kappa\widehat{f}(\omega) $ (see \cite[Theorem 2.9]{MeerschaertsabzikarSPA}). Thus, we can extend the definition of tempered fractional derivatives to a suitable class of functions in $L^2(\R)$ in a natural way, as described below.  For any $\kappa>0$ (and $\lambda>0$), define the fractional Sobolev space
\begin{equation}\label{def:fracSovolev}
W^{\kappa,2}(\R):= \Big\{f\in L^2(\R):\int_{\R} (\lambda^2+\omega^2)^{\kappa}|\hat f(\omega)|^2\,d\omega<\infty \Big\} ,
\end{equation}
which is a Banach space with norm $\|f\|_{\kappa,\lambda}=\|(\lambda^2+\omega^2)^{\kappa/2}\hat f(\omega)\|_2$.  The space $W^{\kappa,2}(\R)$ is the same for any $\lambda>0$ (typically, we take $\lambda=1$)  and all the norms $\|f\|_{\kappa,\lambda}$ are equivalent, since $1+\omega^2\leq\lambda^2+\omega^2\leq \lambda^2(1+\omega^2)$ for all $\lambda\geq 1$, and $\lambda^2+\omega^2\leq1+\omega^2\leq \lambda^{-2}(1+\omega^2)$ for all $0<\lambda<1$.
\begin{defn}{\rm \label{TFDdef2}
The positive (respectively, negative) \textit{tempered fractional derivative} ${\mathbb D}^{\kappa,\lambda}_{\pm}f(t)$ of a function $f\in W^{\kappa,2}(\R)$ is defined as the unique element of $L^2(\R)$  with Fourier transform $\widehat{f}{(\omega)}(\lambda\pm i\omega)^{\kappa}$ for any $\kappa>0$ and any $\lambda>0$.
}\end{defn}

Tempered fractional integrals or derivatives are useful in developing stochastic analysis based on TFLP and TFLP $I\!I$, since we can naturally reexpress these processes based on the former. In fact, for $t<0$, let ${\bf 1}_{[0,t]}(y) := - {\bf 1}_{[-t,0]}(y), y\in \rr$. As shown in Lemma \ref{lem:TFLP connection with fractional calculus}, for $d > 0$ and $t \in \R$, we can write
\begin{equation}\label{e:TFLP-temp_frac_integ_d>0}
S^{I}_{d,\lambda}(t) = \int_{-\infty}^{\infty}\Big({\mathbb I}^{d,\lambda}_{-}{\bf 1}_{[0,t]}-\lambda{\mathbb I}^{d+1,\lambda}_{-}{\bf 1}_{[0,t]}\Big)(y)\ dL(y)
\end{equation}
and
\begin{equation}\label{e:TFLPII-temp_frac_integ_d>0}
S^{I\!I}_{d,\lambda}(t)=\int_{-\infty}^{\infty}\Big({\mathbb I}^{d,\lambda}_{-}{\bf 1}_{[0,t]}\Big)(y)dL(y).
\end{equation}
Likewise, for $-\frac{1}{2}<d<0$ and $t \in \R$,
\begin{equation}\label{e:TFLP-temp_frac_integ_-1/2<d<0}
S^{I}_{d,\lambda}(t) = \int_{-\infty}^{\infty} \Big( {\mathbb D}^{-d,\lambda}_{-}{\bf 1}_{[0,t]}(y)-\lambda
{\mathbb I}^{d+1,\lambda}_{-}{\bf 1}_{[0,t]}(y)\Big) \ dL(y)
\end{equation}
and
\begin{equation}\label{e:TFLPII-temp_frac_integ_-1/2<d<0}
S^{I\!I}_{d,\lambda}(t)=\int_{-\infty}^{\infty}\Big({\mathbb D}^{-d,\lambda}_{-}{\bf 1}_{[0,t]}\Big)(y)dL(y).
\end{equation}
In light of expressions \eqref{e:TFLP-temp_frac_integ_d>0}--\eqref{e:TFLPII-temp_frac_integ_-1/2<d<0}, we are now in a position to construct the theory of stochastic integration with respect to TFLP $I\!I$. Recall that we focus on integration with respect to TFLP \textit{II} because the claims for TFLP are analogous to those for TFBM (see Remark \ref{r:stoch_integ_wrt_TFLP}). Let
\begin{equation}\label{eq:elementarydefn}
f(u)=\sum^{n}_{i=1}a_{i}{\1_{[t_{i},t_{i+1})}(u)}
\end{equation}
be a \textit{step}, or \textit{elementary}, function, where $\{a_i\}_{i=1,\hdots,n}$, $\{t_i\}_{i=1,\hdots,n+1}$, are real numbers such that $t_i \leq a_i \leq t_{i+1}$ for any $i$. Also, let ${\mathcal E}$ be the space of step functions. It is natural to define the stochastic integral of $f \in {\mathcal E}$ with respect to $S^{I\! I}_{d,\lambda}$ by means of the Riemann-Stieltjes-like expression
\begin{equation}\label{eq:THPintegraldefn}
{\Im}^{d,\lambda}(f)=\int_{\mathbb{R}}f(x)dS^{I\!I}_{d,\lambda}(x)=\sum^{n}_{i=1}a_i \Big[ S^{I\!I}_{d,\lambda}(t_{i+1})-S^{I\!I}_{d,\lambda}(t_{i}) \Big].
\end{equation}
Therefore, ${\Im}^{d,\lambda}(f)$ is an infinitely divisible random variable with mean zero.

We first consider the memory parameter range $d  > 0$. It follows immediately from \eqref{e:TFLPII-temp_frac_integ_d>0} that we can write
\begin{equation*}
{\Im}^{d,\lambda}(f)= \int_\rr
\Big({\mathbb{I}}^{d,\lambda}_{-}f\Big)(x)\ dL(x).
\end{equation*}
Moreover, the isometry \eqref{eq:isometry} implies that, for any $f,g\in {\mathcal E}$,
\begin{equation}\label{eq:innerproductTFI}
\begin{split}
&\ip{{\Im}^{d,\lambda}(f)}{{\Im}^{d,\lambda}(g)}_{L^2(\Omega)} =\mathbb{E}\left(\int_{\mathbb{R}}f(x)dS^{I\!I}_{d,\lambda}(x)\int_{\mathbb{R}}g(x)dS^{I\!I}_{d,\lambda}(x)\right) \\
   &= \E L(1)^2  \int_\rr
\Big({\mathbb{I}}^{d,\lambda}_{-}f\Big)(x)
\Big({\mathbb{I}}^{d,\lambda}_{-}g\Big)(x)\ dx .
\end{split}
\end{equation}
In view of expression \eqref{eq:innerproductTFI}, we define and characterize the class of integrands ${ \mathcal{A} }_{1}$ as follows.

\begin{thm}\label{thm:stochasticcalculus for TFI and moving general}
Given $d > 0$ (and $\lambda > 0$), let
\begin{equation}\label{eq:A1 star class}
{\mathcal{A}}_{1}=\Big\{f\in L^{2}(\mathbb{R}):\int_{\mathbb{R}}
\left|\Big({\mathbb{I}}^{d,\lambda}_{-}f\Big)(x)\right|^{2} dx < \infty \Big\}.
\end{equation}
Then, the class of functions ${\mathcal{A}}_{1}$ is a linear space with inner product
\begin{equation}\label{eq:productTFIf}
\begin{split}
{\langle f,g \rangle}_{{\mathcal{A}}_{1}}
&:={\langle F,G\rangle}_{L^{2}(\mathbb{R})},
\end{split}
\end{equation}
where
\begin{equation}\label{eq:defnF}
F(x)= \Big({\mathbb{I}}^{d,\lambda}_{-}f\Big)(x)\qquad {\rm and}\ G(x)= \Big({\mathbb{I}}^{d,\lambda}_{-}g\Big)(x) .
\end{equation}
The set of elementary functions ${\mathcal{E}}$ is dense in ${\mathcal{A}}_{1}$. Moreover, the linear space ${\mathcal{A}}_{1}$ is not complete.
\end{thm}
Note that, although ${\mathcal{A}}_{1} = L^{2}(\mathbb{R})$, the two spaces are endowed with different inner products.

We now define the stochastic integral with respect to TFLP $I\!I$ for any function in ${\mathcal{A}}_{1}$ in the case where $ d > 0$.
\begin{defn}\label{defn:TFI of general f}
\rm{
For any $ d > 0  $  (and $\lambda > 0$),
\begin{equation}\label{eq:TFIof f resp $THP^{1}$ genral}
{\Im}^{d,\lambda}(f) =\int_{\mathbb{R}}f(x)dS^{I\!I}_{d,\lambda}(x) := \int_\rr \Big( \mathbb{I}^{d,\lambda}_{-}f \Big)(x)\ dL(x), \quad f \in {\mathcal{A}}_{1},
\end{equation}
where ${\mathcal{A}}_{1}$ is given by \eqref{eq:A1 star class}.
}
\end{defn}

\begin{rem}
{\rm If one were instead to use the completion $\overline{\mathcal{A}_1}$ of $\mathcal{A}_1$ as a class of integrands, a random element $X\in\overline{\rm Sp}(S^{I\!I}_{d,\lambda})$ could only be represented up to equivalence classes of sequences in $\overline{\mathcal{A}_1}$.  See \cite{PipirasTaqqu} for a detailed discussion. }
\end{rem}

In the following theorem, we establish the link between integrands and stochastic integrals when $d > 0$.
\begin{thm}\label{thm:SLRDisometric}
For any $ d > 0$ (and $\lambda > 0$), the stochastic integral $\Im^{d,\lambda}$ in \eqref{eq:TFIof f resp $THP^{1}$ genral} is an isometry between  ${\mathcal{A}}_{1}$ and a strict subset of
\begin{equation}\label{Spspace}
\overline{\rm Sp}(S^{I\!I}_{d,\lambda})=\Big\{X \in L^2(\Omega):\| {\Im}^{d,\lambda}(f_n)- X\|_{L^2(\Omega)}\rightarrow 0 \textnormal{ for some sequence $(f_n)_{n \in \N} \subseteq \mathcal{E}$}\Big\}.
\end{equation}
%Since ${\mathcal{A}}_{1}$ is not complete, \BCBGD{The spaces ${\mathcal{A}}_{1}$  and $\overline{\rm Sp}(S^{I\!I}_{d,\lambda})$} are not isometric.
\end{thm}

As a consequence of Theorems \ref{thm:stochasticcalculus for TFI and moving general} and \ref{thm:SLRDisometric}, for the memory parameter range $d > 0$ the stochastic integral \eqref{eq:TFIof f resp $THP^{1}$ genral} is well defined as a $L^2(\Omega)$ limit of stochastic integrals constructed from elementary functions.

We now tackle the memory parameter range $-\frac 12< d < 0$. As usual, we first consider integrands in the space of elementary functions ${\mathcal E}$. It follows from \eqref{e:TFLPII-temp_frac_integ_-1/2<d<0} that the stochastic integral \eqref{eq:THPintegraldefn} can be written in the form
\begin{equation*}
{\Im}^{d,\lambda}(f)
=\int_\rr \Big( \mathbb{D}^{-d,\lambda}_{-}f \Big)(x) dL(x), \quad  f\in {\mathcal E}.
\end{equation*}
%Then, ${\Im}^{d,\lambda}(f)$ is an infinitely divisible random variable with mean zero, such that
Moreover, by the isometry \eqref{eq:isometry},
\begin{equation}\label{eq:innerproductTFD}
\begin{split}
\ip{{\Im}^{d,\lambda}(f)}{{\Im}^{d,\lambda}(g)}_{L^2(\Omega)} &=\mathbb{E}\left(\int_\rr \Big(\mathbb{D}^{-d,\lambda}_{-}f \Big)(x)  dL(x)\int_\rr \Big(\mathbb{D}^{-d,\lambda}_{-}g \Big)(x)  dL(x)\right) \\
   &=\int_\rr \Big( {\mathbb{D}}^{-d,\lambda}_{-}f\Big)(x) \Big( {\mathbb{D}}^{-d,\lambda}_{-}g \Big)(x)\ dx,
\end{split}
\end{equation}
for any $f,g\in {\mathcal E}$. In light of expression \eqref{eq:innerproductTFD}, we define and characterize the class of integrands ${ \mathcal{A} }_{2}$ as follows.
%Expression \eqref{eq:innerproductTFD} suggests the following space of integrands for TFLP $I\!I$ in the case $-\frac 12 < d < 0 $, see Theorem \ref{thm:stochasticcalculus for TFI and TFD and moving general} below.
%Define the fractional Sobolev space as the following:  For any $s>0$, we define fractional Sobolev space
%\begin{equation}\label{fracSovolev}
%W^{s,2}(\rr) : = \Big\{ f\in L^2 (\rr) : \int_\rr (\lambda^2 + \omega^2)^s |\hat{f} (\omega)|^2 d\omega <\infty \Big\}
%\end{equation}
%which is a Banach Space with norm $ \|f\|_{s,\lambda} = \| (\lambda^2 + \omega^2)^{s/2} \hat{f}(\omega) \|_2 $.
\begin{thm}\label{thm:stochasticcalculus for TFI and TFD and moving general}
For any $-\frac 12 < d < 0$ (and $\lambda > 0$), let
\begin{equation}\label{eq:A2class}
{\mathcal{A}}_{2} =\Big\{f\in W^{-d,2}(\rr):\varphi_f={\mathbb{D}}^{-d,\lambda}_{-}f\ {\rm for\ some}\ \varphi_f\in{L}^{2}(\rr ) \Big\} .
\end{equation}
Then, the class of functions ${\mathcal A}_2$ is a linear space with inner product
\begin{equation}\label{eq:productTFDf}
{\langle f,g \rangle}_{{\mathcal{A}}_{2}}
:= {\langle {\varphi_f,\varphi_g} \rangle}_{L^{2}(\rr)}.
\end{equation}
The set of elementary functions ${\mathcal{E}}$ is dense in ${\mathcal{A}}_{2}$.  Moreover, the linear space ${\mathcal{A}}_{2}$ is complete.
\end{thm}
%\begin{proof}
%Equation \eqref{eq:EconvPT} is proven in \cite[Lemma 5.1]{PipirasTaqqu}.  For any $L>0$, that proof constructs a sequence of elementary functions $f_n$ such that $\hat f_n(k)\to \1_{[-1,1]}(k)$ almost everywhere on $-L\leq x\leq L$, and shows that $|\hat f_n(k)|\leq C\min\{1,|k|^{-1}\}$ for all $k\in\R$ and all $n\geq 1$. In the notation of that paper, we have $\hat f_n(k)=k^{-1}U_n(k)$.   Apply the dominated convergence theorem to see that
%\[\int_{-L}^{+L} |\hat f_n(k)- \1_{[-1,1]}(k)|^2dk\to 0\]
%and note that
%\[\int_{|k|>L} |\hat f_n(k)- \1_{[-1,1]}(k)|^2 dk\leq 2C^2\int_L^\infty \frac{dk}{k^2}\leq \frac{2C^2}{L} .\]
%Since $L$ is arbitrary,  it follows that $\hat f_n(k)\to \1_{[-1,1]}(k)$ in $L^2(\R)$, and then the result follows as in \cite[Lemma 5.1]{PipirasTaqqu}.
%\end{proof}

We now define the stochastic integral with respect to TFLP $I\!I$ for any function in ${\mathcal{A}}_{2}$ in the case where $-\frac 12 < d < 0$.

\begin{defn}\label{defn:TFD of general f}
{\rm
For any $-\frac 12< d < 0$ (and $\lambda > 0$),
\begin{equation}\label{eq:TFDof f resp TFBM genral}
\Im^{d,\lambda}(f) = \int_\rr f(x)dS^{I\!I}_{d,\lambda}(x):=  \int_\rr \Big( \mathbb{D}^{-d,\lambda}_{-}f \Big)(x)\ dL(x), \quad f\in{\mathcal{A}}_{2},
\end{equation}
where ${\mathcal{A}}_{2}$ is given by \eqref{eq:A2class}.
}
\end{defn}

In the following theorem, we establish the link between integrands and stochastic integrals when $-1/2 < d < 0$. In contrast with the range $d > 0$ (see Theorem \ref{thm:SLRDisometric}), in this case every element in the space $\overline{\rm Sp}(S^{I\!I}_{d,\lambda})$ can be represented as a stochastic integral of a single integrand function $f \in {\mathcal A}_2$.
\begin{thm}\label{thm:SLRDisometric2}
For any $-\frac 12 < d < 0$ (and $\lambda > 0$), the space ${\mathcal{A}}_{2}$ is isometric to $\overline{\rm Sp}(S^{I\!I}_{d,\lambda})$, where $\overline{\textnormal{Sp}}(S^{I\!I}_{d,\lambda})$ is given by \eqref{Spspace}.
%\begin{equation}
%\overline{Sp}(S^{I\!I}_{d,\lambda})=\left\{X:{\Im}^{d,\lambda}(f_n)\rightarrow X\quad in\quad{\mathbb{L}}^{2}(\mathbb{R}),\quad\textrm{for a sequence}
%\left\{f_n\right\}\subset{\mathcal{C}}\right\}.
%\end{equation}
\end{thm}

Note that an element $X\in\overline{\rm Sp}(S^{I\!I}_{d,\lambda})$ is an infinitely divisible random variable. In fact, the law of $S^{I\!I}_{d,\lambda}$ is the limit of infinitely divisible laws and, hence, likewise for $X$. In addition, it has mean zero and finite variance
\begin{equation*}
\Var(X)=\lim_{n\to\infty} \Var[{\Im}^{d,\lambda}(f_n)]
\end{equation*}
(cf.\ \cite{RajputRosinski}, Theorem 2.7).
Moreover, $X$ can be associated with an equivalence class of sequences of elementary functions $(f_n)_{n \in \N}$ such that $\|{\Im}^{\alpha,\lambda}(f_n)- X\|_{L^{2}(\Omega)} \rightarrow 0$ as $n \rightarrow \infty$. Theorem \ref{thm:SLRDisometric2}
states that for any $X\in\overline{\rm Sp}(S^{I\!I}_{d,\lambda})$, there exists a unique $f\in L^2(\R)$ such that $\|f_n - f \|_{L^2(\R)} \rightarrow 0$ as $n \rightarrow \infty$, and that we can write $X=\int_{\mathbb{R}} f(x) dS^{I\!I}_{d,\lambda}(x)$.

%\begin{rem}
%{\rm Theorem \ref{thm:SLRDisometric2} shows that the space ${\mathcal{A}}_{2}$ is isometric to $\overline{\rm Sp}(S^{I\!I}_{d,\lambda})$ when $-\frac 12 < d <0 $. This is in contrast with the case $0 < d < \frac 12 $; see Theorem \ref{thm:SLRDisometric}.   }
%\end{rem}

\begin{rem}\label{r:stoch_integ_wrt_TFLP}
{\rm Stochastic integration with respect to TFLP leads to properties that are analogous to those contained in Theorems 3.5, 3.9, 3.10 and 3.14 in \cite{MeerschaertsabzikarSPA} for the Gaussian case (TFBM). Moreover, these properties can be established by adapting the second order arguments used in \cite{MeerschaertsabzikarSPA}. For the reader's convenience, we summarize the main statements, where $L = \{L(t)\}_{t \in \R}$ is given by \eqref{e:L(t)}.

Let
\begin{equation}\label{e:-1/2<d<0,lambda>0}
-1/2 < d < 0, \quad \lambda > 0.
\end{equation}
Then, the class of functions
$$
{\mathcal A}_3 := \Big\{f \in L^2(\R): \int_{\R} \Big| \I^{-d,\lambda}_{-}f(x) - \lambda \I^{-d+1,\lambda}_{-}f(x)\Big|^2 < \infty\Big\}
$$
is a linear space with inner product $\langle f,g\rangle_{{\mathcal A}_3}:= \langle F,G\rangle_{L^2(\R)}$, where
$$
F(x) = \Gamma(-d + 1) \Big[\I^{-d,\lambda}_{-}f(x) - \lambda \I^{-d+1,\lambda}_{-}f(x)\Big], \quad G(x) = \Gamma(-d + 1) \Big[\I^{-d,\lambda}_{-}g(x) - \lambda \I^{-d+1,\lambda}_{-}g(x)\Big].
$$
Moreover, the space ${\mathcal A}_3$ is not complete. Under \eqref{e:-1/2<d<0,lambda>0}, we define
\begin{equation}\label{e:stoch_integ_FLP_-1/2<d<0}
\int_{\R}f(x)S^{I}_{d,\lambda}(dx):= \Gamma(-d + 1)  \int_{\R} \Big[ \I^{-d,\lambda}_{-}f(x) - \lambda \I^{-d+1,\lambda}_{-}f(x)\Big]dL(x), \quad f \in {\mathcal A}_3.
\end{equation}
Then, the stochastic integral in \eqref{e:stoch_integ_FLP_-1/2<d<0} is an isometry from ${\mathcal A}_3$ into $\overline{\textnormal{Sp}}(S^{I}_{d,\lambda})$. Since ${\mathcal A}_3$ is not complete, these two spaces are not isometric.

Now let
\begin{equation}\label{e:0<d<1/2,lambda>0}
0 < d < 1/2, \quad \lambda > 0,
\end{equation}
and consider the fractional Sobolev space $W^{d,2}(\R)$ as given by \eqref{def:fracSovolev}. Then, the class of functions
$$
{\mathcal A}_4 := \Big\{f \in W^{d,2}(\R): \varphi_f = \D^{d,\lambda}_f - \lambda \I^{1-d,\lambda}_{-}f \textnormal{ for some }\varphi_f \in L^2(\R) \Big\}
$$
is a linear space with inner product $\langle f,g\rangle_{{\mathcal A}_3}:= \langle F,G\rangle_{L^2(\R)}$, where
$$
F(x) = \Gamma(1 -d) \Big[\D^{d,\lambda}_{-}f(x) - \lambda \I^{1-d,\lambda}_{-}f(x)\Big], \quad G(x) = \Gamma( 1-d) \Big[\D^{d,\lambda}_{-}g(x) - \lambda \I^{1-d,\lambda}_{-}g(x)\Big].
$$
Moreover, the space ${\mathcal A}_4$ is not complete. Under \eqref{e:0<d<1/2,lambda>0}, we define
\begin{equation}\label{e:stoch_integ_FLP_0<d<1/2}
\int_{\R}f(x)S^{I}_{d,\lambda}(dx):= \Gamma(1-d)  \int_{\R} \Big[ \D^{d,\lambda}_{-}f(x) - \lambda \I^{1-d,\lambda}_{-}f(x)\Big]dL(x), \quad f \in {\mathcal A}_4.
\end{equation}
Then, the stochastic integral in \eqref{e:stoch_integ_FLP_0<d<1/2} is an isometry from ${\mathcal A}_4$ into $\overline{\textnormal{Sp}}(S^{I}_{d,\lambda})$. Since ${\mathcal A}_4$ is not complete, these two spaces are not isometric.}
\end{rem}

\section{Conclusion}\label{s:conclusion}

In this work, we use exponential tempering to construct two flexible parametric classes of second order, non-Gaussian transient anomalous diffusion models called TFLP and TFLP $I\!I$. In particular, their increment processes exhibit semi-long range dependence, namely, their autocovariance functions decay hyperbolically over small lags and exponentially fast over large lags. We establish the covariance and sample path regularity properties of the TFLP and TFLP $I\!I$ classes. Moreover, with the purpose of constructing a stochastic analysis framework, we use tempered fractional derivatives and integrals to develop the theory of stochastic integration with respect to TFLP and TFLP $I\!I$, which may not be semimartingales.

The results in this paper open up several new research directions. The developed theory provides mathematical tools for the study of solutions of TFLP and TFLP $I\!I$-driven Langevin-type equations. Moreover, it can also be applied in constructing functional limit theorems for unit root problems (cf.\ \cite{FQP}). %In another direction, recall that a stochastic process $X(t)$ is said to be a volatility modulated L\'{e}vy-driven Volterra process if it can be written as $X(t) = \int_0^t g(t,s)\sigma(s)dL(s)$, where $g$ is a deterministic function and $\sigma$ is a predictable stochastic process. In \cite{Barndorff-Nielsen}, stochastic integration with respect to $X(t)$ is developed by identifying the class of appropriate integrands $Y$ such that $\int_0^t Y(s) dX(s)$ is well-defined. The construction of stochastic integration theory with respect to $X(t)$ when $g = g^{I}_{d,\lambda}$ or $g^{I\!I}_{d,\lambda}$ remains an open problem.
From a modeling standpoint, it remains as a future research topic to develop efficient inferential methods for the analysis of geophysical flow and nanobiophysical data. A related research direction is that of the assessment and development of new simulation methods for the TFLP families. This is especially important for TFLP $I\!I$, since the additional integral term in the kernel $g^{I\!I}_{d,\lambda,t}$ makes Stieltjes-based simulation rather computationally costly.

\setcounter{secnumdepth}{0} \section{Acknowledgments}
Farzad Sabzikar would like to thank Alex Lindner for fruitful discussions leading to some results of the paper.
Gustavo Didier was partially supported by the prime award no.\ W911NF--14--1--0475 from the Biomathematics subdivision of the Army Research Office, USA.

%\begin{center}
%\includegraphics[width=.5\linewidth]{sample_paths_alph_16_lam_005.eps}
%\end{center}

%\appendix
%% resets counters for the appendix
\setcounter{section}{0}
\setcounter{secnumdepth}{1}
\renewcommand\thesection{\Alph{section}}
\section{Proofs}\label{Sec5}%\setsectioncounter{A}

\noindent {\it Proof of Proposition \ref{prop:covariance_TFLP}:} The proof of \eqref{eq:TFLPacvf} follows by a similar argument of Proposition 2.3 in \cite{Meerschaertsabzikar} and hence we omit the details.
To show \eqref{e:limitbehaviorofTFLP}, apply the covariance function formula \eqref{eq:TFLPacvf} in Proposition \ref{prop:covariance_TFLP} for $s=t$ to arrive at
\begin{equation}\label{eq:varianceTFLP1}
{\rm Var}\big[S^{I}_{d,\lambda}(t) \big] = \frac{ \E( L(1)^2 ) }{ \Gamma(1+d)^2 } \Bigg[ \frac{ 2\Gamma(1+2d) }{ (2\lambda)^{1+2d} } - \frac{2\Gamma(1+d)}{ \sqrt{\pi} } \Big( \frac{1}{ 2\lambda }\Big)^{d+\frac{1}{2} } |t|^{d+\frac{1}{2}} K_{d+\frac{1}{2}}   (\lambda t)\Bigg].
\end{equation}
The second term inside the bracket tends to zero as $t\to\infty$, since
$$ K_{d+\frac{1}{2}}   (\lambda t) \sim \sqrt{\frac{\pi}{2\lambda t}} e^{-\lambda t}.
$$
Hence, relation \eqref{e:limitbehaviorofTFLP} holds, as claimed. \hfill $\Box$\\
%\begin{equation}
%\lim {\rm Var} \big[ S^{I}_{d,\lambda} (t) \big] = \frac{2 \E( L(1)^2 ) \Gamma(1+2d)}{ \Gamma(1+d)^2 (2\lambda)^{1+2d} }.
%\end{equation}
%This completes the proof.   \hfill $\Box$

\noindent {\it Proof of Proposition \ref{p:modification_Theorem_TFLP1}:} Starting from the definition of TFLP, we can use integration by parts (see \cite{Marquardt}, p.\ 1106) to write
\begin{equation}\label{modification_TFLP_PART1}
\begin{split}
\Gamma(d+1) S^{I}_{d,\lambda}(t) &=  \int_{\rr}\big[ e^{-\lambda(t-x)_{+}} (t-x)_{+}^{d} - e^{-\lambda(-x)_{+}} (-x)_{+}^{d} \big]\ dL(x)\\
&= \int_{-\infty}^{t-} e^{-\lambda(t-x)} (t-x)^{d} dL(x) - \int_{-\infty}^{0-} e^{\lambda x} (-x)^{d} dL(x) \\
&= \lim_{u\uparrow t}\Big( e^{-\lambda(t-u)} (t-u)^{d} L(u) - \int_{-\infty}^{u} L(u) d(e^{-\lambda(t-u)} (t-u)^{d})  \Big) \\
&- \lim_{u\uparrow 0}\Big( e^{\lambda u} (-u)^{d} L(u) - \int_{-\infty}^{u} L(u) d(e^{\lambda u} (-u)^{d})  \Big).
\end{split}
\end{equation}

Using \cite[Proposition 47.11]{Sato}, we have $e^{\lambda v}L(v)\to 0$ as $v \to 0$. Hence, for $d>0$,
\begin{equation*}
\lim_{u\uparrow t} e^{-\lambda(t-u)} (t-u)^{d} L(u) = \lim_{u\uparrow 0} e^{\lambda u} (-u)^{d} L(u) = 0.
\end{equation*}
%Therefore, we can continue \eqref{modification_TFLP_PART1} as follows:
Therefore, we can reexpress \eqref{modification_TFLP_PART1} as
\begin{equation*}
\begin{split}
%\Gamma(d+1) S^{I}_{d,\lambda}(t) &= - \lim_{u\uparrow t} \int_{-\infty}^{u} L(u) \Big(-d e^{-\lambda(t-u)} (t-u)^{d-1} + \lambda e^{-\lambda(t-u)} (t-u)^{d} \Big) du\\
& - \lim_{u\uparrow t} \int_{-\infty}^{u} L(u) \Big(-d e^{-\lambda(t-u)} (t-u)^{d-1} + \lambda e^{-\lambda(t-u)} (t-u)^{d} \Big) du\\
&+ \lim_{u\uparrow 0} \int_{-\infty}^{u} L(u) \Big(-d e^{\lambda u} (-u)^{d-1} + \lambda e^{\lambda u} (-u)^{d} \Big) du\\
&= d \int_{\rr} L(u) \big[ e^{-\lambda(t-u)_{+}} (t-u)_{+}^{d-1} - e^{-\lambda(-u)_{+}} (-u)_{+}^{d-1} \big]\ du\\
&-\lambda \int_{\rr} L(u) \big[ e^{-\lambda(t-u)_{+}} (t-u)_{+}^{d} - e^{-\lambda(-u)_{+}} (-u)_{+}^{d} \big]\ du.\\
\end{split}
\end{equation*}
Hence, \eqref{e:SI_improper_Riemann} holds.

To show the continuity of the process \eqref{e:SI_improper_Riemann}, without loss of generality fix $t\in (a,b) \subseteq \R_+$.  Rewrite the first term in the expression \eqref{e:SI_improper_Riemann} as
$$
\Big\{\int^{a}_{-\infty} +\int^{t}_a \Big\}\Big( e^{-\lambda(t-x)_{+}}{(t-x)_{+}^{d-1}}-e^{-\lambda(-x)_{+}}{(-x)_{+}^{d-1}}\Big)\ L(x)\ dx.
$$
We want to show that this expression is continuous as a function of $t$. On one hand, the mapping $t\mapsto\int_{-\infty}^a L(u) \big[ e^{-\lambda(t-u)_{+}} (t-u)_{+}^{d-1} - e^{-\lambda(-u)_{+}} (-u)_{+}^{d-1} \big]\ du$ is continuous. This is a consequence of the dominated convergence theorem, since
 $$
 \mathbf{1}_{(-\infty,a]}(u)| L(u)|\left| \big[ e^{-\lambda(t-u)_{+}} (t-u)_{+}^{d-1} - e^{-\lambda(-u)_{+}} (-u)_{+}^{d-1}\right|
 $$
 $$
 \leq \mathbf{1}_{(-\infty,a]}(u)| L(u)| \Big( e^{-\lambda(a-u)} (b-u)_{+}^{d-1}+ e^{-\lambda(-u)_{+}} (-u)_{+}^{d-1}\Big)\in L^1(\R),
 $$
 where we use the fact that $L$ is locally bounded. On the other hand, by making the change of variable $z = t-u$,
 \begin{equation}\label{e:int[a,t]e*t^(d-1)dt}
\int^t_a L(u)  e^{-\lambda(t-u)_{+}} (t-u)_{+}^{d-1} \ du = \int_{\R} 1_{[0,t-a]}(z) L(t-z) e^{-\lambda z} z^{d-1}  \ dz.
\end{equation}
However, the integrand in \eqref{e:int[a,t]e*t^(d-1)dt} is bounded in absolute value by
$$
\sup_{w \in (a,b)}|L(w)| 1_{[0,b-a]}(z) e^{-\lambda z} z^{d-1} \in L^1(\R).
$$
Therefore, by the dominated convergence theorem, the mapping $ t \mapsto \int_a^t L(u) \big[ e^{-\lambda(t-u)_{+}} (t-u)_{+}^{d-1} - e^{-\lambda(-u)_{+}} (-u)_{+}^{d-1} \big]\ du$ is also continuous. Hence, the first term in the expression \eqref{e:SI_improper_Riemann} is continuous as a function of $t$, as claimed. Again by the dominated convergence theorem, the second term in the expression \eqref{e:SI_improper_Riemann} is also continuous as a function of $t$. This establishes that the process \eqref{e:SI_improper_Riemann} is continuous. \hfill $\Box$\\
\noindent {\it Proof of Theorem \ref{t:TFLP_sample_path}:}  %\BCBcomment{Please see the remark at the end of the proof.}
First, we establish (a). We use the modification of $S^{I}_{d,\lambda}$ given in Theorem \ref{p:modification_Theorem_TFLP1} to write
\begin{equation}\label{prop29generalparta}
\begin{split}
|S^{I}_{d,\lambda}(t) - S^{I}_{d,\lambda}(s)|&\leq \frac{1}{\Gamma(d)} \int_{\rr} \Big| e^{-\lambda(t-u)_{+}} (t-u)_{+}^{d-1} - e^{-\lambda(s-u)_{+}} (s-u)_{+}^{d-1} \Big| \Big|L(u)\Big| \ du\\
%&\leq \frac{1}{\Gamma(d)} \int_{\rr} \Big| e^{-\lambda(t-u)_{+}} (t-u)_{+}^{d-1} - e^{-\lambda(s-u)_{+}} (s-u)_{+}^{d-1} \Big| \Big|L(u)\Big| \ du\\
&+ \frac{\lambda}{\Gamma(d+1)} \int_{\rr} \Big| e^{-\lambda(t-u)_{+}} (t-u)_{+}^{d} - e^{-\lambda(s-u)_{+}} (s-u)_{+}^{d} \Big| \Big|L(u)\Big| \ du.
\end{split}
\end{equation}
Recall that $0 < d \leq 1/2$. For notational simplicity, consider a parameter $\beta$, which can be interpreted either as $d$ or $d-1$, i.e. $\beta \in (-1,-1/2]\cup (0,1/2]$.
 Define
\begin{equation*}
W_{\beta}(s,t) = \int_{\rr} \Big| e^{-\lambda(t-u)_{+}} (t-u)_{+}^{\beta} - e^{-\lambda(s-u)_{+}} (s-u)_{+}^{\beta} \Big| \Big|L(u)\Big| \ du.
\end{equation*}
%Since $-T \leq s \leq t \leq T$, we obtain
For $s$, $t$ satisfying  $-T \leq s \leq t \leq T$, we obtain
\begin{equation*}
\begin{split}
W_{\beta}(s,t) &= \int_{s}^{t} e^{-\lambda(t-u)} (t-u)^{\beta}  \Big|L(u)\Big| \ du
 + \int_{-\infty}^{s} \Big| e^{-\lambda(t-u)} (t-u)^{\beta} - e^{-\lambda(s-u)} (s-u)^{\beta} \Big| \Big|L(u)\Big| \ du\\
& \leq \sup_{|u|\leq T} |L(u)| \int_{s}^{t} e^{-\lambda(t-u)} (t-u)^{\beta}  \ du + \int_{-\infty}^{s}  e^{-\lambda(t-u)} \Big| (t-u)^{\beta} - (s-u)^{\beta} \Big| \Big|L(u)\Big| \ du\\
&+ \int_{-\infty}^{s}  (s-u)^{\beta} \Big| e^{-\lambda(t-u)} - e^{-\lambda(s-u)}  \Big| \Big|L(u)\Big| \ du.
\end{split}
\end{equation*}
Using the substitution $h=t-s$, we get
\begin{equation}\label{prop29generalpart}
\begin{split}
W_{\beta}(s,t) & \leq \frac{h^{\beta+1}}{\beta+1}\sup_{|u|\leq T} |L(u)| + e^{-\lambda h}\int_{-\infty}^{s}  e^{-\lambda(s-u)} \Big| (h+s-u)^{\beta} - (s-u)^{\beta} \Big| \Big|L(u)\Big| \ du\\
&+ \int_{-\infty}^{s} | e^{-\lambda h}-1 | (s-u)^{\beta} e^{-\lambda(s-u)}  \Big|L(u)\Big| \ du\\
&= \frac{h^{\beta+1}}{\beta+1}\sup_{|u|\leq T} |L(u)| + e^{-\lambda h}\int_{0}^{\infty}  e^{-\lambda v} \Big| (h+v)^{\beta} - v^{\beta} \Big| \Big|L(s-v)\Big| \ dv\\
&+ ( 1- e^{-\lambda h} ) \int_{0}^{\infty}  v^{\beta} e^{-\lambda v}  \Big| L(s-v) \Big| \ dv\\
&=: I_1 + I_2 + I_3.
\end{split}
\end{equation}
%\BCBcommentnew{It appears there was a typo in the equation above, there was a term $e^{-\lambda(1-u)}$ that seemed to be incorrect.  The change appears in blue above.}
Since $L$ is locally bounded, %\GDcomment{As a consequence of what condition?} \BCBcommentnew{It's from Doob's maximal inequality.}
then
\begin{equation}\label{prop29I1part}
I_1 \leq C_{1}(\omega) h^{\beta+1}
\end{equation}
for an almost surely finite random variable $C_1$. Next, observe that
\begin{equation}\label{Sato_prop}
\limsup_{|v|\to\infty}\hspace{0.5mm}\frac{|L(v)|}{|v|} = 0
%\lim\sup_{|v|\to\infty} \frac{|L(v)|}{|v|} = 0
\end{equation}
by \cite[Proposition 48.9]{Sato}.
In particular, the integrands appearing in $I_2$ and $I_3$ are finite almost surely (since $\lambda>0$). Since $( 1- e^{-\lambda h} ) \leq \lambda h$ for $h>0$, %\BCBcommentnew{I removed the $h\in[0,2\pi]$ here since it seemed unnecessary.}
we conclude that there is an almost sure finite continuous random variable $C_3(\omega)$ such that
\begin{equation}\label{prop29I3part}
I_{3}(\omega) \leq C_{3}(\omega) h
\end{equation}
 for all $-T\leq s\leq t\leq T$.

In regard to $I_2$, consider the decomposition
\begin{equation}\label{e:I2=int01+int1infty}
\Big\{\int_{0}^{1} + \int_{1}^{\infty} \Big\} \hspace{1mm}e^{-\lambda v} | (h+v)^{\beta} - v^{\beta} | |L(s-v)| \ dv.
\end{equation}
%we will use the mean value theorem, for each $v>0$ there exists some $v_{h}\in [v,v+h]$ such that $(h+v)^{\beta} - v^{\beta} = h\beta
By the mean value theorem, for each $v>0$ there exists some $v_{h}\in [v,v+h]$ such that $(h+v)^{\beta} - v^{\beta} = h\beta
{v_{h}}^{\beta-1}$. Thus, we can bound the second integral in \eqref{e:I2=int01+int1infty} by
\begin{equation}\label{prop29I21part}
\begin{split}
&\int_{1}^{\infty}  e^{-\lambda v} | (h+v)^{\beta} - v^{\beta} | |L(s-v)| \ dv \\
&\leq \int_{1}^{\infty}  e^{-\lambda v} h \beta \max\{v^{\beta}, (v+h)^{\beta}\} |L(s-v)| \ dv \leq C_{2,1} h
%&\leq C_{2,1} h_{1}
\end{split}
\end{equation}
$-T \leq s\leq t \leq T$. %where  $C_{2,1}$ is an almost surely finite random variable by \eqref{Sato_prop}. On the other hand,
In \eqref{prop29I21part}, $C_{2,1}$ is an almost surely finite random variable as a consequence of \eqref{Sato_prop}. On the other hand, the first integral in \eqref{e:I2=int01+int1infty} can be bounded by
\begin{equation*}
\begin{split}
&\int_{0}^{1}  e^{-\lambda v} \Big| (h+v)^{\beta} - v^{\beta} \Big| \Big|L(s-v)\Big| \ dv \\
&\leq \sup_{v\in[-T-1,T]} |L(v)| \int_{0}^{1} \Big| (h+v)^{\beta} - v^{\beta} \Big|  \ dv \\
&= \sup_{v\in[-T-1,T]} |L(v)|  \Big| \int_{0}^{1} (h+v)^{\beta}\ dv - \int_0^1 v^{\beta}   \ dv \Big| \\
&= \sup_{v\in[-T-1,T]} |L(v)|\hspace{0.5mm} \frac{1}{\beta+1} \Big| (1+h)^{\beta+1} -  h^{\beta+1} -1 \Big| \\
& \leq \sup_{v\in[-T-1,T]} |L(v)| \hspace{0.5mm}\frac{1}{\beta+1} \Big( |(1+h)^{\beta+1}-1| +  h^{\beta+1} \Big).
\end{split}
\end{equation*}
Using a Taylor expansion, it follows that there is an almost surely finite random variable $C_{2,2}$ such that
\begin{equation}\label{prop29I22part}
\begin{split}
&\int_{0}^{1}  e^{-\lambda v} \Big| (h+v)^{\beta} - v^{\beta} \Big| \Big|L(s-v)\Big| \ dv \leq   C_{2,2}|h|^{\min (1,\beta+1)}
%&\leq   C_{2,2}|h|^{\min (1,\beta+1)}
\end{split}
\end{equation}
for $s,t\in [-T,T]$. Combining \eqref{prop29generalpart}, \eqref{prop29I1part}, \eqref{prop29I3part}, \eqref{prop29I21part}, and \eqref{prop29I22part},
%together
we see that %\BCBcommentnew{If we change the parameter range for $\beta$ here the ``min'' should be removed below and in the previous equation.  Also, it seems that this proof could be cleaned up by first using the stationarity of the increments, possibly.  Also: we actually seem to get $d$-holder continuity, instead of the current weaker statement( $\gamma$-Holder for every $\gamma\in(0,d)$), no?}
\begin{equation}\label{e:Walpha(s,t)=<Ch^min(1,alpha+1)}
|W_{\beta}(s,t)| \leq C_{\beta} h^{\min(1,\beta+1)}
\end{equation}
for $s,t\in [-T,T]$, where $C_{\beta}$ is an almost surely finite random variable. Applying \eqref{e:Walpha(s,t)=<Ch^min(1,alpha+1)} to \eqref{prop29generalparta} with $\beta = d$ and $\beta = d-1$ yields
%then give
\begin{equation*}
%|S^{I}_{d,\lambda}(t) - S^{I}_{d,\lambda}(s)| \leq C_{T} |t-s|^{\min(1,d)} \qquad {\rm for }\quad s,t\in[-T,T]
|S^{I}_{d,\lambda}(t) - S^{I}_{d,\lambda}(s)| \leq C_{T} |t-s|^{d},  \quad s,t\in[-T,T],
\end{equation*}
%which is desired result.
which establishes \eqref{TFLPsamplepath}.% \BCBcomment{I don't understand the above proof.  I believe the proof should be a direct application of Kolmogorov-$\breve{\textnormal{C}}$entsov (see proof below):}\\

To show (b), let $-\frac{1}{2} < d < 0$. In this case, the kernel function $g^{I}_{d,\lambda,\cdot}(s)$ is not locally bounded and in fact the mapping $t\longmapsto g^{I}_{d,\lambda,t}(s)$, $t\in\R$, is unbounded and discontinuous for all $s$. Therefore, Theorem 4 in \cite{Rosinski89} implies that the sample paths of $S^{I}_{d,\lambda}$ are unbounded and discontinuous with positive probability, as claimed. \hfill $\Box$ \\%and this proves (b).

\noindent {\it Proof of Proposition \ref{prop:TFLN}:} To prove (a), note that TFLN has the same covariance structure as tempered fractional Gaussian noise (TFGN), up to a constant. Expression \eqref{e:gammaI(h)_decay_semi-LRD} can be obtained by following the same argument as in Chen et al.\ \cite[Appendix 2]{chen:wang:deng:2017} for the asymptotic behavior of TFGN over large covariance lags.

To show (b), let $a(t)$ be the time domain kernel of the moving average representation \eqref{eq:TFGNmoving} of TFLN. Then, the spectral density is given by
$$
h^{I}(\omega) = \frac{1}{2\pi} \hspace{1mm} \Big| \int_{\R}e^{-i\omega t} a(t) dt\Big|^2  = \frac{1}{2\pi} \hspace{1mm}\Big| \frac{ e^{i\omega} -1 }{ (\lambda + i\omega)^{d+1}  }   \Big|^2.
$$
This establishes \eqref{e:TFLN_density}. \hfill $\Box$\\

%%
%Using Proposition \ref{prop:covariance_TFLP} and \eqref{eq:TFLNdef}, for $h>0$, one can verify that
%\begin{equation}\label{eq:asymptotic TFLN1}
%\begin{split}
%\gamma(h) &= -\frac{ \E[L(1)^2 ]}{\Gamma(d+1) (2\lambda)^{d+\frac 12} }\\
%&\bigg[ (h+1)^{d+\frac 12} K_{ d+\frac 12 }(\lambda (h+1)) + (h-1)^{d+\frac 12} K_{d+\frac 12}(\lambda (h-1))
%-2 h^{ d+\frac 12 } K_{d+\frac 12}(\lambda h) \bigg]
%\end{split}
%\end{equation}
%It is known that $K_{\nu}(x) \sim \frac{e^{-x}\sqrt{\pi}}{\sqrt{2x}}$ as $x\to\infty$ and then
%\begin{equation}\label{eq:asymptotic TFLN2}
%\begin{split}
%&(h+1)^{d+\frac 12} K_{ d+\frac 12 }(\lambda (h+1)) + (h-1)^{d+\frac 12} K_{d+\frac 12}(\lambda (h-1))
%-2 h^{ d+\frac 12 } K_{d+\frac 12}(\lambda h)\\
%&= h^{d+\frac 12} \bigg[ (1+\frac 1h)^{d+\frac 12} K_{ d+\frac 12 }(\lambda (h+1)) + (1-\frac 1h)^{d+\frac 12} K_{d+\frac 12}(\lambda (h-1))
%-2 K_{d+\frac 12}(\lambda h) \bigg]\\
%& \sim h^{d+\frac 12} \bigg[ \frac{ \sqrt{\pi} e^{-\lambda h} h^{-\frac 12} }{ \sqrt{2\lambda} } ( e^{\lambda} + e^{-\lambda} - 2 ) \bigg]
%\end{split}
%\end{equation}
%as $h\to \infty$. Now combining \eqref{eq:asymptotic TFLN1} and \eqref{eq:asymptotic TFLN2} together yields
%\begin{equation*}
%\gamma(h) \sim C(d,\lambda) h^{d} e^{-\lambda h},
%\end{equation*}
%where $C(d,\lambda) = -\frac{ \E[L(1)^2 ]}{\Gamma(d+1) (2\lambda)^{(d+1)} }(e^{-\lambda} + e^{\lambda} - 2 ) \sim -\frac{\lambda^2 \E[L(1)^2 ]}{\Gamma(d+1) (2\lambda)^{(d+1)} }$
%and this completes the proof.
%
%

The next lemma is mentioned in Section \ref{s:TFLPII}. As a consequence of the lemma, $S^{I\!I}_{d,\lambda}(t)$ is well defined for any $t>0$.
\begin{lem}\label{lem:g squar integrable}
Let $g^{I\! I}_{d,\lambda,t}(y)$ be the function \eqref{hdef0}. Then,
\begin{equation}\label{e:g squar integrable}
g^{I\! I}_{d,\lambda,t}(y) \in L^2(\R)
\end{equation}
for any $t \in \R $ and any $\lambda >0$, $d > -\frac{1}{2}$.
\end{lem}
\noindent {\it Proof of Lemma \ref{lem:g squar integrable}: } {Let $t >0$. By applying Minkowski's inequality to \eqref{hdef0}, we arrive at
$$
\|  g^{I\! I}_{d,\lambda,t}(\cdot)\|_{2}
\le \Big(\int_\R (t-y)_+^{2d} \ e^{-2 \lambda (t-y)_+} \ dy\Big)^{1/2} + \Big(\int_\R (-y)_+^{2d} \ e^{-2 \lambda (-y)_+} \ dy\Big)^{1/2}
$$
$$
+ \lambda  \Big( \int_\R \ \Big\{\int^t_0 (s-y)_+^{d} \ e^{-\lambda (s-y)_+} \ ds \Big\}^{2}dy \Big)^{1/2}  < \infty,
$$
where finiteness is a consequence of the facts that $2d+1>0$ and $\lambda >0$. Since $g^{I\! I}_{d,\lambda,-t}(y) = -g^{I\! I}_{d,\lambda,t}(y+t) $ for any $t, y \in \R$, \eqref{e:g squar integrable} holds.\hfill $\Box$\\

%%%%%%%%%%%%%%%%%%%%%%%%%%%%%%%%%%
%{\it Proof of Lemma \ref{lem:g squar integrable}: } Let $t >0$ and apply Minkowski's inequality for \eqref{hdef0} to see
%$\|  h_{d,\lambda}(t;\cdot)\|^2_{2}
%\le \int (t-y)_+^{2d} \ e^{-2 \lambda (t-y)_+} \ dy + \int (-y)_+^{2d} \ e^{-2 \lambda (-y)_+} \ dy
%+ \lambda^2  \Big\{\int_0^t \ ds \big(\int (s-y)_+^{2d} \ e^{-2 \lambda (s-y)_+} \ dy \big)^{1/2} \Big\}^2 $
%$< \infty $ due to $ 2d+1 >0 $ and $\lambda >0$.
%Since $h_{d,\lambda}(-t;y) = - h_{d,\lambda}(t;y+t) $ for any $t, y \in \R$ this proves
%the Lemma. \hfill $\Box$
%%%%%%%%%%%%%%%%%%%%%%%%%%%%%%%%%

\noindent {\it Proof of Proposition \ref{prop:THPcovariance}:} We first note that $g^{I\! I}_{d,\lambda,t}(y) = d\int_{0}^{t} (s-y)_+^{d-1} e^{-\lambda (s-y)_+} ds $, where $g^{I\! I}_{d,\lambda,t}(y)$ is the function given by \eqref{hdef0}. Hence,
\begin{equation}\label{eq:TFLPI_second_rep}
S^{I\!I}_{d,\lambda}(t)=\frac{1}{\Gamma(d+1)} \int_{\rr} g^{I\! I}_{d,\lambda,t}(y)\ dL(y)
= \frac{1}{\Gamma(d)} \int_{\rr}\int_{0}^{t} (s-y)_+^{d-1} e^{-\lambda (s-y)_+} ds \ dL(y)
\end{equation}
%It is known that, \BCBGD{for $f \in L^2(\R)$,}
From Proposition \ref{prop:RajputRosinski},
\begin{equation}\label{eq:covarianceTaqqu}
{\text {Cov}}\Big(\int_{\rr}f(y)\ dL(y), \int_{\rr}g(y)\ dL(y)\Big)=\E[L(1)^2] \int_{\rr} f(y)g(y) dy
\end{equation}
%\BCBGD{(see, e.g., \cite{Marquardt}, p.\ 1103).}
Now, by Lemma \ref{lem:g squar integrable}, we can apply \eqref{eq:covarianceTaqqu} to TFLP $I\!I$ in \eqref{eq:TFLPI_second_rep} to write
\begin{equation}\label{eq:COVARIANCE1}
\begin{split}
&{\text {Cov}}\Big(S^{I\!I}_{d,\lambda}(t),S^{I\!I}_{d,\lambda}(s)\Big)=\frac{\E[L(1)^2]}{(\Gamma(d))^{2}}\int_{\rr} g^{I\! I}_{d,\lambda,t}(y)g^{I\! I}_{d,\lambda,s}(y) dy\\
&=\frac{\E[L(1)^2]}{(\Gamma(d))^{2}}\int_{\rr}\Big(\int_{0}^{t}\int_{0}^{s}(u-y)_{+}^{d-1}(v-y)_{+}^{d-1}e^{-\lambda(u-y)_{+}}e^{-\lambda(v-y)_{+}}dv\ du\Big)dy\\
&=\frac{\E[L(1)^2]}{(\Gamma(d))^{2}}\int_{0}^{t}\int_{0}^{s}\Bigg[\int_{-\infty}^{\min(u,v)}
(u-y)^{d-1}(v-y)^{d-1}e^{-\lambda (u-y)}e^{-\lambda (v-y)}dy\Bigg]dv\ du.
\end{split}
\end{equation}
Using the relation
\begin{equation}\label{eq:Bessel}
\int_{0}^{\infty}x^{\nu-1}(x+\beta)^{\nu-1}e^{-\mu x}dx=\frac{1}{\sqrt{\pi}}\left(\frac{\beta}{\mu}\right)^{\nu-\frac{1}{2}}e^{\frac{\beta \mu}{2}}\ \Gamma(\nu)K_{\frac{1}{2}-\nu}\Big(\frac{\beta\mu}{2}\Big),
\end{equation}
(see \cite{Gradshteyn}, p.\ 348),
%\eqref{eq:Bessel},
%\begin{equation}\label{eq:standard integral formula}
%\int_{0}^{\infty}x^{\nu-1}(x+\beta)^{\nu-1}e^{-\mu x}\ dx=\frac{1}{\sqrt{\pi}}\big(\frac{\beta}{\mu}\big)^{\nu-\frac{1}{2}}e^{\frac{\beta\mu}{2}}
%\Gamma(\nu)K_{\frac{1}{2}-\nu}(\frac{\beta\mu}{2}),
%\end{equation}
%for $|\arg \beta|<\pi$, Re $\mu>0$, Re $\nu>0$,
we have
\begin{equation}\label{eq:COVARIANCE2}
%\begin{split}
\int_{-\infty}^{\min(u,v)}(u-y)^{d-1}(v-y)^{d-1}e^{-\lambda (u-y)}e^{-\lambda (v-y)}dy
%&=\int_{-\infty}^{u}(u-y)^{d-1}(v-y)^{d-1}\ e^{-\lambda (u-y)}e^{-\lambda (v-y)}\1_{\{y<u<v\}}dy\\
%&+\int_{-\infty}^{v}(u-y)^{d-1}(v-y)^{d-1}\ e^{-\lambda (u-y)}e^{-\lambda (v-y)}\1_{\{y<v<u\}}dy\\
%&=\int_{0}^{\infty}{w}^{d-1}(w+v-u)^{d-1}\ e^{-\lambda w}e^{-\lambda (w+v-u)}\1_{\{w<u<v\}}dw\\
%&+\int_{0}^{\infty}{w}^{d-1}(w+u-v)^{d-1}\ e^{-\lambda w}e^{-\lambda (w+u-v)}\1_{\{w<v<u\}}dw\\
%&=\int_{0}^{\infty}{w}^{d-1}(w+|u-v|)^{d-1}\ e^{-2\lambda w}\ e^{-\lambda |u-v|}dw\\
=\frac{\Gamma(d)}{\sqrt{\pi}}\Big(\frac{|u-v|}{2\lambda}\Big)^{d-\frac{1}{2}}K_{d-\frac{1}{2}}(\lambda|u-v|).
%\end{split}
\end{equation}
Therefore, from \eqref{eq:COVARIANCE1} and \eqref{eq:COVARIANCE2}, we have
\begin{equation*}
{\text {Cov}}\Big(S^{I\!I}_{d,\lambda}(t),S^{I\!I}_{d,\lambda}(s)\Big)=\frac{\E[L(1)^2]}{\sqrt{\pi}\Gamma(d)(2\lambda)^{d-\frac{1}{2}}}\int_{0}^{t}\int_{0}^{s}
|u-v|^{d-\frac{1}{2}}K_{d-\frac{1}{2}}(\lambda|u-v|)dv\ du
\end{equation*}
for any $d>0$ and $\lambda>0$, as claimed.  \hfill $\Box$\\

\noindent {\it Proof of Proposition \ref{modification_Theorem_TFLP2}:} The proof follows the similar technique that was employed in Theorem \ref{p:modification_Theorem_TFLP1} and hence we omit it.   \hfill $\Box$\\
%First, we note that $S^{I\!I}_{d,\lambda}$ given by \eqref{eq:defTFLP second} can be written as
%\begin{equation*}
%\Gamma(d)S^{I\!I}_{d,\lambda}(t) = \int_{-\infty}^{t}\int_{0}^{t} (s-x)_+^{d-1} e^{-\lambda(s-x)_+} ds\ dL(x)
%\end{equation*}
%provided $d>0$. Now, using integration by parts (see \cite{Marquardt}, p.\ 1106) on the aforementioned equation yields
%\begin{equation*}
%\begin{split}
%\Gamma(d)S^{I\!I}_{d,\lambda}(t) & = \Bigg( L(x) \int_{0}^{t} (s-x)_+^{d-1} e^{-\lambda(s-x)_+} ds \Bigg)_{-\infty}^{t^-}\\
%&-  \int_{-\infty}^{t} L(x)\int_0^t\Big[ (1-d) (s-x)_+^{d-2} e^{-\lambda(s-x)_+} + \lambda (s-x)_+^{d-1} e^{-\lambda(s-x)_+}\Big]ds dx \\
%\end{split}
%\end{equation*}
%\begin{equation*}
%= (d-1) \int_{-\infty}^{t} \int_0^t (s-x)_+^{d-2} e^{-\lambda(s-x)_+} dsL(x) dx
%- \lambda \int_{-\infty}^{t}\int_0^t (s-x)_+^{d-1} e^{-\lambda(s-x)_+}ds L(x) dx.
%\end{equation*}
%Continuity follows similarly as in Theorem \ref{p:modification_Theorem_TFLP1}. Hence, the claim holds. \hfill $\Box$\\

\noindent {\it Proof of Theorem \ref{t:TFLPI bounds}:} We use the Kolmogorov-$\breve{\textnormal{C}}$entsov theorem (e.g., \cite{karatzas:shreve:2000}, p.\ 53) to establish the claim. Since $\lambda >0$ is fixed, we can assume $\lambda =1 $ without loss of generality. Since the increments of $S^{I\!I}_{d,1}(t)$ are stationary, it suffices to show that
\begin{equation}\label{e:E|SIId,1(t)|2=<Ct(1+beta)}
\E |S^{I\!I}_{d,1}(t)|^2  \le C t^{1+\beta}
\end{equation}
for some $ \beta > 0$ and all $ 0 < t < 1 $. Consider $g^{I\! I}_{d,1,t}$ as in \eqref{eq:defTFLP second}. By \eqref{eq:covarianceTaqqu},
\begin{equation*}
\E |S^{I\!I}_{d,1}(t)|^2  = C \int_{-\infty}^t (g^{I\! I}_{d,1,t}(y))^2  \   dy =:  C (I_1 + I_2),
\end{equation*}
where

\begin{equation*}
I_1 = \int_{-t}^t (g^{I\! I}_{d,1,t}(y))^2  \ dy= \frac{1}{\Gamma(d)}\int_{-t}^t\left( \int_{0}^{t} (s-x)^{d-1} e^{-(s-x)} ds\right)^2dx$$$$\le C\int_{-t}^{t} (t - y)^{2d} \ dy \le C t^{2d+1}
\end{equation*}

and
\begin{equation*}
\begin{split}
I_2 = \int_{-\infty}^{-t} (g^{I\! I}_{d,1,t}(y))^2  \ dy &\le C\int_t^\infty ((t + y)^{d} \ e^{-t-y}
- y^{d} \ e^{-y})^2  \ dy\\
&\  +  C\int_t^\infty \Big\{\int_0^t (s+y)^{d} \ e^{-s-y} \ ds \Big\}^2 \ dy
=  C(I_2' + I_2'').
\end{split}
\end{equation*}
Using  $|(t + y)^{d} \ e^{-(t-y)} - y^{d} \ e^{-y}| \le |\ e^{-t}-1|\, \ e^{-y} (t + y)^{d} +  \ e^{-y}\, |(t + y)^{d}
- y^{d} | \le C t \, \ e^{-y} (t + y)^{d} +  C t \, \ e^{-y} y^{d} $ we obtain
$I'_2 \le C t^2 $ and, similarly, $I''_2 \le C t^2 $, implying
$I_1 + I_2 \le C(t^{2d+1} + t^2) $ and
$ \E |S^{I\!I}_{d,1}(t)|^2  \le C(t^{2d+1} + t^2)\leq t^{2d+1}$ since $d\in(0,1/2]$ and $0 < t < 1$. Hence, \eqref{e:E|SIId,1(t)|2=<Ct(1+beta)} is satisfied with $\beta = 2d$. This completes the proof.
To show (b), note that when $-\frac{1}{2} < d < 0$ $g^{I\! I}_{d,\lambda,\cdot}(s)$ is not locally bounded and $t\longmapsto g^{I\! I}_{d,\lambda,t}(s)$, $t\in\R$ is unbounded and discontinuous for all $s$, and so the same proof in part (b) of Theorem \ref{t:TFLP_sample_path} applies.
\hfill $\Box$\\

\noindent {\it Proof of Proposition \ref{p:TFLNII}:}
%%%%%%%%%%%%%%%%%
%\noi (a) For a TFGNII, Sabzikat and Surgailis \cite{SurgailisFarzadTFSMII} showed that $\gamma(h) \asymp  e^{-\lambda h } h^{ d -1 }$
%as $h\to\infty$. From \eqref{eq:isometry}, one can see TFBMII and TFLPI have the same second order structure up to some constants and this gives the desired result in (a).
%%%%%%%%%%%%%%%%%
To show (a), note that the autocovariance function of a TFGN $I\! I$ satisfies $\gamma(h) \asymp  e^{-\lambda h } h^{ d -1 }$ as $h\to\infty$ (see \cite{SurgailisFarzadTFSMII}). From \eqref{eq:isometry}, TFBM $I\! I$ and TFLP $I\! I$ have the same second order structure up to constants. Hence, \eqref{e:gamma(h)_semi-LRD_TFLNII} holds.

%\noi(b)It is known that the sum of a covariance stationary can be obtained from the value of the spectral density at the origin. That is, $\sum_{h\in \Z}\gamma(h) = 2\pi f(0)$ where
%$f(\omega)$ is the spectral density of the process. Sabzikat and Surgailis \cite[Proposition 3.1]{SurgailisFarzadTFSMII} obtained the spectral density of TFGNII and showed that
%$f(0) = \lambda^{2d} / 2\pi$. Again, since TFBMII and TFLPI have the same second order structure up to some constants, then we can conclude that the increments of these processes, TFGNII and TFLNI, have the same
%behavior. Therefore $\sum_{h\in \Z}\gamma(h) = 2\pi f(0) = \lambda^{2d}$ for $ d>0 $ and this shows that TFLNI has short memory which completes the proof of part (b).
%%%%%%%%%%%%%%%%%%%%
%For part(b), use the fact that
%$$X^{I\!I}(t) = \int_{\rr} \Big( \mathbb{I}^{d,\lambda}_{-}{\bf 1}_{[t,t+1]}(x) \Big) \ dL(x), $$ where $\Big(\mathbb{I}^{d,\lambda}_{-} f\Big)(x)$ is TFI of function $f$. Now, the spectral density follows by:
%\begin{equation*}
%\begin{split}
%h^{I\!I}(\omega) &= \frac{\sigma^2}{2\pi} \Big|  e^{-i\omega t} \Big( \mathbb{I}^{d,\lambda}_{-}{\bf 1}_{[t,t+1]}(\omega) \Big) dt   \Big|^2\\
%&= \frac{\sigma^2}{2\pi} \Big| (\lambda + i\omega)^{-d}\int_{t}^{t+1} e^{-i\omega x} dx \Big|^2
%= \frac{\sigma^2}{2\pi} \frac{ 2( 1- \cos(\omega) ) }{(\lambda^2 + \omega^2)^{d}\ \omega^2}
%\end{split}
%\end{equation*}
%and this completes the proof of part (b) and the proposition. \hfill $\Box$
%%%%%%%%%%%%%%%%%%%%%%%%%%%

To show (b), let $\Big(\mathbb{I}^{d,\lambda}_{-} f\Big)(x)$ be as in \eqref{eq:positivetempered fractional integral} with $\kappa = d$. Note that the process $X^{I\!I}_{d,\lambda}$ as in \eqref{eq:TFGNmoving} has the integral representation
$$
X^{I\!I}_{d,\lambda}(t) = \int_{\rr} \Big( \mathbb{I}^{d,\lambda}_{-}{\bf 1}_{[t,t+1]}(x) \Big) \ dL(x).
$$
Therefore, its spectral density is given by
\begin{equation*}
\begin{split}
h^{I\!I}(\omega) &= \frac{1}{2\pi} \Big|\int_{\R}  e^{-i\omega t} \Big( \mathbb{I}^{d,\lambda}_{-}{\bf 1}_{[t,t+1]}(\omega) \Big) dt   \Big|^2\\
&= \frac{1}{2\pi} \Big| (\lambda + i\omega)^{-d}\int_{t}^{t+1} e^{-i\omega x} dx \Big|^2
= \frac{1}{2\pi} \frac{ 2( 1- \cos(\omega) ) }{(\lambda^2 + \omega^2)^{d}\ \omega^2},
\end{split}
\end{equation*}
as claimed. \hfill $\Box$\\
\noindent{{\it Proof of Proposition \ref{p:not_a_semimg}}: Write
\begin{equation}\label{e:phiI_phiII}
\phi^I(x)= x_+^de^{-\lambda x_+}, \quad \phi^{II}(x)=   x_+^d e^{-\lambda x_+} + \lambda \int_0^{x} u_+^d e^{-\lambda u_+} du,
\end{equation}  and note that
$$
S^I_{d,\lambda}(t) = \frac{1}{\Gamma(d+1)}\int_{\R} \{\phi^I(t-x)-\phi^I(-x)\}dL(x),
$$
\begin{equation}\label{e:S1,S2,phiI,phiII}
S^{I\! I}_{d,\lambda}(t) = \frac{1}{\Gamma(d+1)}\int_{\R} \{\phi^{II}(t-x)-\phi^{II}(-x)\}dL(x).
\end{equation}
For $x>0$, the derivatives $\eta^I(x):=\frac{d}{dx}\phi^I(x)$, $\eta^{II}(x):=\frac{d}{dx}\phi^{II}(x)$ exist and satisfy $\eta^I(x)\sim \eta^{II}(x)\sim d x^{d-1}$, $x \to 0^+$.  Hence
$$
\int_a^b \left|\eta^I(x)\right|^\alpha dx =\infty,\quad  \int_a^b \left|\eta^{II}(x)\right|^\alpha dx=\infty
$$
for any interval $[a,b)$ containing 0 whenever $\alpha(d-1)+1<0$, i.e., whenever $d+\frac{1}{\alpha}<1$. Hence, by Corollary 3.4 in \cite{BassePedersen},  the processes
$$
\int_0^t\phi^I(t-x)dL(x), \quad \int_0^t\phi^{II}(t-x)dL(x), \quad t\geq 0
$$
are not $(\mathcal{F}^{L}_t)_{t\geq0}$-semimartingales, where $(\mathcal{F}^{L}_t)_{t\geq0}=\sigma\{L(s); 0\leq s \leq t \}$.  Thus, since $L$ is symmetric, in view of the representations \eqref{e:S1,S2,phiI,phiII}, by Lemma 5.2 of \cite{BassePedersen} $S^I_{d,\lambda}$ and $S^{I\! I}_{d,\lambda}$ are not  $(\mathcal{F}^{L,\infty}_t)_{t\geq 0}$-semimartingales.} \hfill $\Box$
\\
%\BCBcomment{we could shorten the proof below, too: ``The proof is similar to that of Theorem 3.9 in \cite{cheridito:2004}; by an application of the Stochastic Fubini theorem, one arrives at... (expression), as was to be shown." }
\noindent {\it Proof of Proposition \ref{semimartingaleTFLP}}:  The proof is similar to that of Theorem 3.9 in \cite{cheridito:2004}.  Write $\eta^I(x)=\frac{d}{dx}\phi^I(x)$ where $\phi^I$ is given in \eqref{e:phiI_phiII}.  Note since $d>1/2$, $\eta^I\in L^2(\R)$, and hence the integral $\int_\R \eta^I(x)  dL(x)$ is well-defined.   Now,
$$\Gamma(d+1)S^I_{d,\lambda}(t) = \int_{-\infty}^t \{\phi^I(t-x) - \phi^I(-x) \}dL(x)$$
$$
=  \int_{-\infty}^0 \{\phi^I(t-x) - \phi^I(-x)\} dL(x) +  \int_{0}^t \phi^I(t-x)  dL(x)
$$
$$
= \int_{-\infty}^0 \int_0^t \eta^I(s-x) ds dL(x) +  \int_{0}^t \int_{x}^t\eta^I(s-x)ds  dL(x).
$$
Hence, by a stochastic version of the Fubini theorem (e.g. \cite{protter:2003}, Theorem 65), the above process has a version that is equal to
$$
 \int_0^t\int_{-\infty}^0 \eta^I(s-x) dL(x)ds +  \int_{0}^t \int_{0}^s\eta^I(s-x)  dL(x)ds
= \int_0^t\int_{-\infty}^s \eta^I(s-x) dL(x)ds.
$$
This establishes $(i)$.

We now turn to $(ii)$. First note that
\begin{equation*}
S^{I\! I}_{d,\lambda}(t) = \frac{1}{\Gamma(d+1)}\int_{\R}\{ \phi^{II}(t-x)-\phi^{II}(-x)\}dL(x),
%S^{I\!I}_{d,\lambda}(t) = \int_{\rr} \phi(t-y) - \phi(-y) dL(y),
\end{equation*}
where $\phi^{II}$ is given in \eqref{e:phiI_phiII}.  Since  $d>1/2$, $\frac{d}{dx}\phi^{II}(x)\in L^2(\R)$, and the rest of the proof can be done similarly to that of part $(i)$.% Hence, all the claims hold.
\hfill $\Box$\\

%Now, the rest of the proof is similar to that of Proposition \ref{semimartingaleTFLP} and this completes the proof. \hfill $\Box$\\

The following lemma is used in Section \ref{Sec3}.
\begin{lem}\label{lem:TFLP connection with fractional calculus}%with $d>-\frac{1}{2}$ and $\lambda>0$.
Let $ S^{I}_{d,\lambda} $ and $S^{I\!I}_{d,\lambda}(t) $ be a TFLP and TFLP $I\!I$ given by \eqref{eq:TFLPdef} and \eqref{eq:defTFLP second}, respectively.  Then, for every $t\in\R$,
\begin{itemize}
\item[(a)]  when $d>0$, expressions \eqref{e:TFLP-temp_frac_integ_d>0} and \eqref{e:TFLPII-temp_frac_integ_d>0} hold;
\item[(b)] when $-\frac{1}{2}<d<0$, expressions \eqref{e:TFLP-temp_frac_integ_-1/2<d<0} and \eqref{e:TFLPII-temp_frac_integ_-1/2<d<0} hold.
\end{itemize}
\end{lem}
\noindent {\it Proof of Lemma \ref{lem:TFLP connection with fractional calculus}:} The proofs can be developed along the same lines of that of Lemma 3.4 in \cite{MeerschaertsabzikarSPA} for TFLP, and of Proposition 2.5 in \cite{SurgailisFarzadTFSMII} for TFLP $I\!I$. \hfill $\Box$\\

\noindent {\it Proof of Theorem \ref{thm:stochasticcalculus for TFI and moving general} :} To show that ${\mathcal{A}}_{1}$ is an inner product space, it suffices to establish that ${\langle f,f \rangle}_{{\mathcal{A}}_{1}}=0$ implies $f=0$ $dx$--a.e.  If ${\langle f,f \rangle}_{{\mathcal{A}}_{1}}=0$, then in view of \eqref{eq:productTFIf} and \eqref{eq:defnF} we have ${\langle F,F\rangle}_{2}=0$, so $F(x)=\Big({\mathbb{I}}^{d,\lambda}_{-}f\Big)(x)=0$ $dx$--a.e. Then,
\begin{equation}\label{eq1:stochasticcalculus for TFI and moving general}
\Big({\mathbb{I}}^{d,\lambda}_{-}f\Big)(x)=0\quad dx\text{--a.e.}
%\Big({\mathbb{I}}^{d,\lambda}_{-}f\Big)(x)=0\quad\text{for almost every $x\in\R$.}
\end{equation}
Apply ${\mathbb{D}}^{d,\lambda}_{-}$ to both sides of equation \eqref{eq1:stochasticcalculus for TFI and moving general} and use Lemma 2.14 in \cite{MeerschaertsabzikarSPA} to get $f(x)=0$ $dx$--a.e. Hence, ${\mathcal{A}}_{1}$ is an inner product space, as claimed.

Next, we want to show that the set of elementary functions ${\mathcal E}$ is dense in ${\mathcal{A}}_{1} \subseteq L^2(\R)$. For any $f\in{\mathcal{A}}_{1}$, we also have $f\in{L}^{2}(\mathbb{R})$, and hence there exists a sequence of elementary functions $(f_n)_{n \in \N}$ in $L^2(\R)$ such that $\|f-f_n\|_2\to 0$ as $n \rightarrow \infty$. However,
\begin{equation*}
\|f-f_n\|^2_{{\mathcal{A}}_{1}} ={\langle f-f_n,f-f_n \rangle}_{{\mathcal{A}}_{1}}={\langle F-F_n,F-F_n\rangle}_2= \|F-F_n\|^2_2,
\end{equation*}
where
$
F_n(x)=\Big({\mathbb{I}}^{d,\lambda}_{-}{f_n}\Big)(x)
$
and $F(x)$ is given by \eqref{eq:defnF}. It can be further shown that $\|{\mathbb{I}}^{\kappa,\lambda}_{-}(f)\|_{2}\leq C\|f\|_{2}$ for some constant $C$. Then,
\begin{equation*}
\|f-f_n\|_{{\mathcal{A}}_{1}}=\left\|F-F_n\right\|_{2}=\|{\mathbb{I}}^{d,\lambda}_{-}(f-f_n)\|_{2}\leq C\|f-f_n\|_{2}.
\end{equation*}
Since $\|f-f_n\|_2\to 0$ as $n \rightarrow \infty$, it follows that the set of elementary functions is dense in ${{\mathcal{A}}_{1}}$.
Finally, using the example provided in the \cite[Theorem 3.1]{PipirasTaqqu}, one can show that ${\mathcal{A}}_{1}$ is not complete.

The following proposition can be established by a direct adaptation of the proof of Proposition 2.1 in \cite{PipirasTaqqu}.
\begin{prop}\label{p:Prop2.1_PipTaq_for_Levy}
For $d > -1/2$, $\lambda>0$, let ${\mathcal E}$ be the set of elementary functions, let ${\mathcal I}^{d,\lambda}(f)$ be an integral \eqref{eq:TFIof f resp $THP^{1}$ genral} of $f \in {\mathcal E}$ with respect to the L\'{e}vy process $L$ as in \eqref{e:L(t)}. Suppose ${\mathcal D}$ is a set of deterministic functions on $\R$ such that: $(i)$ ${\mathcal D}$ is an inner product space with an inner product $\langle f,g\rangle_{{\mathcal D}}$ for $f,g \in {\mathcal D}$; (ii) ${\mathcal E} \subseteq {\mathcal D}$ and  $\langle f,g\rangle_{{\mathcal D}} = \langle {\mathcal I}^{d,\lambda}(f),{\mathcal I}^{d,\lambda}(g)\rangle_{L^2(\Omega)}$, $f,g \in {\mathcal E}$; $(iii)$ the set is dense in ${\mathcal D}$. Then,
\begin{itemize}
\item [(a)] there is an isometry between the space ${\mathcal D}$ and a linear subspace of $\overline{\textnormal{Sp}}(S^{I\!I}_{d,\lambda} )$ which is an extension of the mapping $f \mapsto {\mathcal I}^{d,\lambda}(f)$, $f \in {\mathcal E}$;
    \item [(b)] ${\mathcal D}$ is isometric to $\overline{\textnormal{Sp}}(S^{I\!I}_{d,\lambda} )$ itself if and only if ${\mathcal D}$ is complete.
\end{itemize}
\end{prop}
%%%%%%%%%%%%%%%%%%%%%%%%%%%%%%%%%%%%%
%Proposition 2.1 in Pipiras and Taqqu \cite{PipirasTaqqu} implies that, if $\mathcal{D}$ is an inner product space such that $(f,g)_{\mathcal{D}}=\ip{{\Im}^{d,\lambda}(f)}{{\Im}^{d,\lambda}(g)}_{L^2(\Omega)}$ for all $f,g\in\mathcal{E}$, and if $\mathcal{E}$ is dense $\mathcal{D}$, then there is an isometry between $\mathcal{D}$  and a linear subspace of $\overline{\rm Sp}( S^{I\!I}_{d,\lambda} )$ that extends the map $f\to{\Im}^{d,\lambda}(f)$ for $f\in\mathcal{E}$, and furthermore, $\mathcal{D}$ is isometric to $\overline{\rm Sp}( S^{I\!I}_{d,\lambda} )$ itself if and only if $\mathcal{D}$ is complete.
%%%%%%%%%%%%%%%%%%%%%%%%%%%%%%%%%%%%%
We are now in a position to prove Theorem \ref{thm:SLRDisometric}.

\noindent {\it Proof of Theorem \ref{thm:SLRDisometric} :} Since $\|{\mathbb{I}}^{\kappa,\lambda}_{-}(f)\|_{2}\leq C\|f\|_{2}$ then the stochastic integral \eqref{eq:TFIof f resp $THP^{1}$ genral} is well-defined for any $f\in{\mathcal{A}}_{1}$. By using the isometry \eqref{eq:isometry} and expression \eqref{eq:TFIof f resp $THP^{1}$ genral}, it follows from Proposition \ref{p:Prop2.1_PipTaq_for_Levy} and \eqref{eq:productTFIf} that, for any $f,g\in{{\mathcal{A}}_{1}}$,
 \[{\langle f,g \rangle}_{{\mathcal{A}}_{1}}={\langle F,G\rangle}_{L^{2}(\mathbb{R})}
 =\ip{{\Im}^{d,\lambda}(f)}{{\Im}^{d,\lambda}(g)}_{L^2(\Omega)}.\]
Then, Theorem \ref{thm:stochasticcalculus for TFI and moving general} implies that ${\mathcal{A}}_{1}$ is isometric to a subset of $\overline{\rm Sp}(S^{I\!I}_{d,\lambda})$, as claimed. However, again by Theorem \ref{thm:stochasticcalculus for TFI and moving general}, ${\mathcal{A}}_{1}$ is not complete. Therefore, ${\mathcal{A}}_{1}$ is isometric to a strict subset of $\overline{\rm Sp}(S^{I\!I}_{d,\lambda})$. \hfill $\Box$\\

Lemmas \ref{prop:sobolevA2} and \ref{lem:EconvPT}, stated and proved next, are used in the proof of Theorem \ref{thm:stochasticcalculus for TFI and TFD and moving general}.

\begin{lem}\label{prop:sobolevA2}
Under the assumptions of Theorem \ref{thm:stochasticcalculus for TFI and TFD and moving general}, every $f\in W^{-d,2}(\R)$ is an element of  ${\mathcal{A}}_{2}$
for $-\frac 12 < d <0$ and $\lambda >0 $, i.e., as sets, $W^{-d,2}(\R)=\mathcal{A}_2$.
\end{lem}
\noindent {\it Proof of Lemma \ref{prop:sobolevA2}:}
Given $f\in W^{-d,2}(\R)$, we need to show that
\begin{equation}\label{eq:SobolevinA2*}
\varphi_f=\mathbb{D}^{-d,\lambda}_{-} f
\end{equation}
for some $\varphi_f\in L^{2}( \rr )$. From the definition \eqref{def:fracSovolev} we see that $\int (\lambda^2+\omega^2)^{-d}|\hat{f}(\omega)|^2\ d\omega<\infty$.  Define  $h_1(\omega)=(\lambda-i\omega)^{-d}\hat{f}(\omega)$ and note that $h_1$ is the Fourier transform of some function $\varphi_1\in L^{2}( \rr )$. Define $\varphi_f:=\varphi_{1}$ so that
\begin{equation}\label{eq:invertfourier}
\widehat{\varphi_f}(\omega) = \widehat{\varphi_1}(\omega) = \widehat{f}(\omega)(\lambda-i\omega)^{-d}.
\end{equation}
Since $f\in W^{-d,2}(\R)\subset L^2(\R)$, we can apply Definition \ref{TFDdef2} to get the desired result. \hfill $\Box$\\

 We state the following lemma that will be used to proof Theorem
\ref{thm:stochasticcalculus for TFI and TFD and moving general}. We refer the reader to \cite[Lemma 3.12]{MeerschaertsabzikarSPA} for the proof of the Lemma.
\begin{lem}\label{lem:EconvPT}
Suppose the assumptions of Theorem \ref{thm:stochasticcalculus for TFI and TFD and moving general} hold. If $f\in W^{-d,2}(\rr)$, then there exists a sequence of functions $(f_n)_{n \in \N} \subseteq {\mathcal E}$ such that $\|f_n - f\|_{L^2(\rr)}$. Moreover, when $-\frac 12 < d <0$,
\begin{equation}\label{eq:EconvPT}
\int_\rr |\widehat{f}_n(\omega)-\widehat f(\omega)|^2|\omega|^{-2d}d\omega \to 0, \quad\text{$n\to\infty$.}
\end{equation}
\end{lem}

\noindent {\it Proof of Theorem \ref{thm:stochasticcalculus for TFI and TFD and moving general} }: For $f\in{\mathcal{A}}_{2}$ we define
\begin{equation}\label{eq:defnnormA2}
\|f\|_{{\mathcal{A}}_{2}}=\sqrt{\ip ff_{{\mathcal{A}}_{2}}}=\sqrt{\ip{\varphi_f}{\varphi_f}_2}=\|{\varphi_f}\|_{2},
\end{equation}
where $\varphi_f$ is given by \eqref{eq:SobolevinA2*}.  Next, use \eqref{eq:invertfourier} to see that
\begin{equation}\label{eq:invertfourier2}
\widehat{\varphi_f}(\omega) = (\lambda-i\omega)^{-d}\widehat{f}(\omega) .
\end{equation}
To verify that \eqref{eq:productTFDf} is an inner product, it suffices to show that, if ${\langle f,f \rangle}_{{\mathcal{A}}_{2}}=0$, then
\begin{equation} \label{e:f=0_dx-a.e.}
f=0 \quad dx\textnormal{--a.e.}
\end{equation}
In fact,
 \begin{equation}\label{eq:fA2*norm}
0 = \|f\|_{{\mathcal{A}}_{2}}^2=\|{\varphi_f}\|_{2}^2=\|{\widehat{\varphi_f}}\|^{2}_{2}
=\int_\rr |\widehat{f}(\omega)|^{2} (\lambda^2+\omega^2)^{-d}\ d\omega
\end{equation}
implies that $\widehat{f}(\omega)=0$ $d\omega$--a.e. Hence, \eqref{e:f=0_dx-a.e.} holds.
%%%%%%%%%%%%%%%%%%%%%%%%%%
%To verify that \eqref{eq:productTFDf} is an inner product, note that if ${\langle f,f \rangle}_{{\mathcal{A}}_{2}}=0$ then
% \begin{equation}\label{eq:fA2*norm}
%\|f\|_{{\mathcal{A}}_{2}}^2=\|{\varphi_f}\|_{2}^2=\|{\widehat{\varphi_f}}\|^{2}_{2}
%=\int_\rr |\widehat{f}(\omega)|^{2} (\lambda^2+\omega^2)^{-d}\ d\omega
%\end{equation}
%equals zero, which implies that $|\widehat{f}(\omega)|=0$ almost everywhere, and then $f=0$ almost everywhere. This proves that \eqref{eq:SobolevinA2*} is an inner product.
%%%%%%%%%%%%%%%%%%%%%%%%%%

We now show that ${\mathcal{E}}$ is dense in ${\mathcal{A}}_2$. By Lemma \ref{lem:EconvPT}, there is a sequence $(f_n)_{n \in \N} \subseteq {\mathcal E }$ such that
\begin{equation}\label{e:|fn-f|L2->0}
\|f_n-f\|_2\to 0, \quad n \rightarrow \infty,
\end{equation}
and \eqref{eq:EconvPT} holds. %It is easy to see that any elementary function is an element of $W^{-d,2}(\rr)$ and then Lemma \ref{prop:sobolevA2} implies that it is also an element of ${\mathcal{A}}_2$.
On the other hand, by Lemma \ref{prop:sobolevA2}, ${\mathcal E} \subseteq W^{-d,2}(\rr) \subseteq {\mathcal{A}}_2$.
By \eqref{eq:fA2*norm}, we can write
\begin{equation*}
\|f_n-f\|^{2}_{{\mathcal{A}}_2} =\int_\rr \Big|\widehat{f}_n(\omega)-\widehat{f}(\omega)\Big|^{2}(\lambda^2 + \omega^2)^{-d}\ d\omega =: I_{1} + I_{2},
\end{equation*}
where
\begin{equation*}
I_1 = \int_{|\omega| < \lambda}\Big|\widehat{f}_n(\omega)-\widehat{f}(\omega)\Big|^{2} (\lambda^2+\omega^2)^{-d}\ d\omega, \quad I_2 = \int_{|\omega| \geq \lambda}\Big|\widehat{f}_n(\omega)-\widehat{f}(\omega)\Big|^{2} (\lambda^2+\omega^2)^{-d}\ d\omega.
\end{equation*}
%\begin{equation}
%I_1 = \int_{|\omega| < \lambda}\Big|\widehat{f}_n(\omega)-\widehat{f}(\omega)\Big|^{2} (\lambda^2+\omega^2)^{-d}\ d\omega
%\end{equation}
%and
%\begin{equation}
%I_2 = \int_{|\omega| \geq \lambda}\Big|\widehat{f}_n(\omega)-\widehat{f}(\omega)\Big|^{2} (\lambda^2+\omega^2)^{-d}\ d\omega.
%\end{equation}
Since $|\omega| < \lambda$, then $I_1  \leq 2\lambda^{-2d} \int_{\rr} |\widehat{f}_n(\omega)-\widehat{f}(\omega)\Big|^{2} \ d\omega\to 0$ as $n\to\infty$, where convergence is a consequence of \eqref{e:|fn-f|L2->0}. Moreover, by \eqref{eq:EconvPT}, $I_2  \leq 2^{-d} \int_{\rr} \Big|\widehat{f}_n(\omega)-\widehat{f}(\omega)\Big| |\omega|^{-2d}\ d\omega\to 0$ as $n\to\infty$.
Hence, $\|f_n-f\|^{2}_{{\mathcal{A}}_2} \to 0$ as $n\to\infty$, namely, ${\mathcal{E}}$ is dense in ${\mathcal{A}}_2$.
%
%
%Next we want to show that ${S}(\mathbb{R})$ is dense in ${{\mathcal{A}}_{2}}$ \fbox{True but how does this help?}.
%\fbox{Dr. Meerschaert, We proved ${S}(\mathbb{R})$ is dense in ${{\mathcal{A}}_{2}}$, }
%\fbox{because we need to show that the set of elementary functions are dense in ${{\mathcal{A}}_{2}}$ }.
%
%
%
%Recall that the fractional Sobolev norm is defined by $\|f\|_{\beta,\lambda}=\|(\lambda^2+k^2)^{\beta/2}\hat f(k)\|_2$.  Then it follows from \eqref{eq:fA2*norm} that
%\begin{equation*}
%\begin{split}
%\|f\|_{{\mathcal{A}}_{2}}&=\int_{-\infty}^{\infty}|\widehat{f}(k)|^{2}\frac{k^2}{(\lambda^2+k^2)^{1-\beta}}\ dk
%\leq \int (\lambda^2+k^2)^{\beta}\left|\widehat{f}(k)\right|^{2}\ dk=\|f\|_{\beta,\lambda}
%\end{split}
%\end{equation*}
%for any $f\in{\mathcal{A}}_{2}$.  Since ${\mathcal S}(\R)$ is dense in $W^{\beta,2}(\R)$ \cite[Proposition 4.11]{Demengel}, for any $f\in{\mathcal{A}}_{2}$ it follows from Lemma \ref{prop:sobolevA2} that there exists a sequence of functions $(f_n)$ in ${\mathcal S}(\R)$ such that $\|f_n-f\|_{{\mathcal{A}}_{2}}\leq \|f_n-f\|_{\beta,\lambda}\to 0$, and so ${\mathcal S}(\R)$ is also dense in the space ${\mathcal{A}}_{2}$.

It only remains to show that ${\mathcal{A}}_{2}$ is complete. In fact, let $\big( f_n \big)_{n \in \N}$ be a Cauchy sequence in ${\mathcal{A}}_{2}$. Then, by using the inner product \eqref{eq:productTFDf}, the corresponding sequence
$\big( { \varphi_{f_n} } \big)_{n \in \N}$ is Cauchy in $L^2(\rr)$. Again by the inner product \eqref{eq:productTFDf}, and since $L^2(\rr)$ is complete, there exists $\varphi_{f^*}$ such that $\|f_n- f^*\|_{{\mathcal{A}}_{2}} = \|{\varphi_{f_n} } - {\varphi_{f^*} }\|_{2}\to 0$, $n\to\infty$. Hence, $f^* \in {\mathcal A}_2$ and ${\mathcal{A}}_{2}$ is complete. \hfill $\Box$\\

%Finally, we want to show that ${\mathcal{A}}_{2}$ is complete. Let $\big( f_n \big)$ be a Cauchy sequence in ${\mathcal{A}}_{2}$. Then using the inner product \eqref{eq:productTFDf},
%the corresponding
%$\big( { \varphi_{f_n} } \big)$ is Cauchy in $L^2(\rr)$ and since $L^2(\rr)$ is complete there exist $\varphi_{f^*}$ such that $\|\varphi_{f_n}-\varphi_{f^*}\|_2\to 0$ as $n\to\infty$. Again, using \eqref{eq:productTFDf} we have
%\begin{equation}\label{eq:defnnormA2}
%\|f_n- f^*\|_{{\mathcal{A}}_{2}} = \|{\varphi_{f_n} } - {\varphi_{f^*} }\|_{2}\to 0,\ {\rm as}\ n\to\infty
%\end{equation}
%and this completes the proof. \hfill $\Box$\\

\noindent {\it Proof of Theorem \ref{thm:SLRDisometric2}:}
By Lemma \ref{prop:sobolevA2}, the stochastic integral \eqref{eq:TFDof f resp TFBM genral} is well-defined for any $f\in{\mathcal{A}}_{2}$.
Since ${\mathcal{A}}_{2}$ is a complete space with inner product \eqref{eq:productTFDf} and $\mathcal{E}$ is dense, then Proposition \ref{p:Prop2.1_PipTaq_for_Levy} implies that ${\mathcal{A}}_{2}$ is isometric to $\overline{\rm Sp}( S^{I\!I}_{d,\lambda} )$. This completes the proof.   \hfill $\Box$

\bibliographystyle{abbrv}
\bibliography{tempered_flp}

\small

\bigskip

\noindent \begin{tabular}{lll}
B. Cooper\ Boniece &  Gustavo Didier &  Farzad Sabzikar \\
Department of Mathematics and Statistics & Mathematics Department & Department of Statistics \\
Washington University in St.\ Louis & Tulane University  & Iowa State University \\
CB 1146, One Brookings Drive & 6823 St.\ Charles Avenue & 2438 Osborn Drive  \\
St. Louis, MO 63130-4899, USA & New Orleans, LA 70118, USA & Ames, IA 50011-1090, USA \\
{\it bcboniece@wustl.edu} & {\it gdidier@tulane.edu} & {\it sabzikar@iastate.edu}  \\
\end{tabular}\\

\smallskip

%\end{thebibliography}

\end{document}